\documentclass[a4paper,11pt,oneside,reqno]{amsart}

\usepackage{hyperref}
\usepackage[headinclude,DIV13]{typearea}
\areaset{15.1cm}{25.0cm}
\usepackage[utf8]{inputenc}
\usepackage[T1]{fontenc}
\usepackage{amsthm,amsmath,txfonts,amssymb,bbm,hyperref,mathtools,csquotes,thmtools,verbatim,xfrac,bigints}
\usepackage{chngcntr}
\usepackage{textcomp}
\usepackage{tikz}
\usepackage{pgfplots}
\usepackage{multicol}

\usepackage[UKenglish]{babel}

\theoremstyle{plain}
\newtheorem{thm}{Theorem}[section]
\newtheorem{cor}[thm]{Corollary}
\newtheorem{prop}[thm]{Proposition}
\newtheorem{lemma}[thm]{Lemma}

\newtheorem{assumption}{Assumption}
\theoremstyle{definition}
\newtheorem{definition}[thm]{Definition}
\theoremstyle{remark}
\newtheorem{rem}[thm]{Remark}

\def\paragraph#1{\noindent \textbf{#1}}

\numberwithin{equation}{section}

\usepackage{graphicx}
\usepackage{xcolor}
\usepackage{color}
\usepackage{tikz}
\usepackage{pgfplots}
\usetikzlibrary{shapes.misc}
\tikzset{cross/.style={cross out, draw=black, minimum size=2*(#1-\pgflinewidth), inner sep=0pt, outer sep=0pt},
	cross/.default={2pt}}
\def\addlegendimage{\csname pgfplots@addlegendimage\endcsname}

\pgfkeys{/pgfplots/number in legend/.style={%
		/pgfplots/legend image code/.code={%
			\node at (0.295,-0.0225){#1};
		},%
	},
}

\definecolor{darkgreen}{rgb}{0,.6,0}
\definecolor{darkagenta}{rgb}{.5,0,.5}
\definecolor{darkred}{rgb}{1,0,0}%was 0.85
\definecolor{darkblue}{rgb}{0,0,.4}
\definecolor{black}{rgb}{0,0,0}
\definecolor{white}{rgb}{1,1,1}

\begin{document}
	\title[{E}xtremal process of the scale-inhomogeneous 2d DGFF]{Extremes of the 2d scale-inhomogeneous discrete {G}aussian free field: {E}xtremal process in the weakly correlated regime}
	\author{Maximilian Fels, Lisa Hartung}
	\address{M. Fels\\Institut f\"ur Angewandte Mathematik\\
		Rheinische Friedrich-Wilhelms-Universität\\ Endenicher Allee 60\\ 53115 Bonn, Germany }
	\email{fels@iam.uni-bonn.de, wt.iam.uni-bonn.de/maximilian-fels}
	\address{L. Hartung\\Institut f\"ur Mathematik\\ Johannes Gutenberg-Universität Mainz\\ Staudingerweg 9\\ 55099 Mainz, Germany}
	\email{lhartung@uni-mainz.de, https://sites.google.com/view/lisahartung/home}
	\thanks{M.F. is funded by the Deutsche Forschungsgemeinschaft (DFG, German Research Foundation) - project-id 211504053 - SFB 1060 and Germany’s Excellence Strategy – GZ 2047/1, project-id 390685813 – “Hausdorff Center for Mathematics” at Bonn University.}
	\maketitle
	\begin{abstract}
		We prove convergence of the full extremal process of the two-dimensional scale-inhomogeneous discrete Gaussian free field in the weak correlation regime. The scale-inhomogeneous discrete Gaussian free field is obtained from the 2d discrete Gaussian free field by modifying the variance through a function $\mathcal{I}:[0,1]\rightarrow [0,1]$. The limiting process is a cluster Cox process. The random intensity of the Cox process depends on the $\mathcal{I}^\prime(0)$ through a random measure $Y$ and on the $\mathcal{I}^\prime(1)$ through a constant $\beta$. We describe the cluster process, which only depends on $\mathcal{I}^\prime(1)$, as points of a standard 2d discrete Gaussian free field conditioned to be unusually high.
	\end{abstract}
	\section{Introduction}
	Log-correlated processes have received a lot of attention in recent years, see e.g. \cite{MR3101852,MR3129797,2015arXiv150304588D,2016arXiv160600510B,MR3164771,MR3361256,fyodorov1,MR3531703,MR3594368,MR3911893}. Prominent examples are branching Brownian motion (BBM), the two-dimensional discrete Gaussian free field (DGFF), cover times of Brownian motion on the torus, characteristic polynomials of random unitary matrices or local maxima of the randomized Riemann zeta function on the critical line.
One of the key features in these models is that their correlations are such that they start to become relevant for the extreme values of the processes. In particular, one is interested in the structure of the extremal processes that arises when the size of the index set tends to infinity. In the case of the 2d DGFF, one considers the field indexed by the vertices of a lattice box of side length $N$, where $N$ is taken to infinity.
In this paper, we study the extremal process of the scale-inhomogeneous 2d DGFF in the weakly correlated regime. The model first appeared as a tool to prove Poisson-Dirichlet statistics of the extreme values of the 2d DGFF \cite{MR3354619}. In the context of the 2d DGFF, it is the natural analogue model of the variable-speed BBM or time-inhomogeneous branching random walk (BRW). We start with a precise definition of the model we consider in the following.
\begin{definition}[2d discrete Gaussian free field (DGFF)]\label{definition:dgff}
	Let $N\in \mathbb{N}$ and $V_N=[0,N)^2\cap \mathbb{Z}^2$. Then, the centred Gaussian field $\{\phi^N_v\}_{v\in V_N}$ with correlations given by the Green kernel
	\begin{align}
		\mathbb{E}\left[\phi^N_v \phi^N_w\right]=G_{V_N}(v,w)\coloneqq\frac{\pi}{2} \mathbb{E}_v\left[\sum_{k=0}^{\tau_{\partial V_N}-1}\mathbbm{1}_{S_k=w} \right],\text{ for } v,w\in V_N
	\end{align}
	is called DGFF on $V_N$. Here, $\mathbb{E}_v$ is the expectation with respect to the SRW $\{S_k\}_{k\geq 0}$ on $\mathbb{Z}^2$ started in $v$ and $\tau_{\partial V_N}$ denotes the stopping time of the SRW hitting the boundary $\partial V_N$.
\end{definition}
\begin{definition}[2d scale-inhomogeneous DGFF]\label{definition:scale_dependent_dgff}
	Let $\{\phi^N_v\}_{v\in V_N}$ be a DGFF on $V_N$. For $v=(v_1,v_2)\in V_N$ and $\lambda\in (0,1)$, set
	\begin{align}
		[v]_\lambda \equiv [v]^N_\lambda \coloneqq \left(\left[v_1-\frac{1}{2}N^{1-\lambda},v_1+\frac{1}{2}N^{1-\lambda}\right]\times \left[v_2-\frac{1}{2}N^{1-\lambda},v_2+\frac{1}{2}N^{1-\lambda}\right] \right)\cap V_N.
	\end{align}
	We set  $[v]^N_0\coloneqq V_N$ and $[v]^N_1 \coloneqq \{v\}.$ We denote by $[v]^o_\lambda$ the interior of $[v]_\lambda$. Let $\mathcal{F}_{\partial{[v]_{\lambda}} \cup [v]_{\lambda}^c } \coloneqq \sigma\left(\{ \phi^N_v, v \notin [v]_{\lambda}^o \}\right)$ be the $\sigma-$algebra generated by the random variables outside $[v]_{\lambda}^o$.
	For $v\in V_N$, let
	\begin{align}\label{equation:condition_dgff__lambda_env}
	\phi^N_v(\lambda)=\mathbb{E}\left[\phi^N_v | \mathcal{F}_{\partial{[v]_{\lambda}} \cup [v]_{\lambda}^c } \right],\quad  \lambda \in [0,1].
	\end{align}
	We denote by $\nabla \phi^N_v(\lambda)$ the derivative $\partial_\lambda \phi^N_v(\lambda)$ of the DGFF at vertex $v$ and scale $\lambda$.
	Moreover, let $s\mapsto \sigma(s)$ be a non-negative function such that $\mathcal{I}_{\sigma^2}(\lambda)\coloneqq\int_{0}^{\lambda}\sigma^2(x)\mathrm{d}x$ is a function on $[0,1]$ with $\mathcal{I}_{\sigma^2}(0)=1$ and $\mathcal{I}_{\sigma^2}(1)=1$.
	The $2$d scale-inhomogeneous DGFF on $V_N$ is a centred Gaussian field, $\psi^N \coloneqq \{\psi^N_v \}_{v \in V_N},$ defined as
	\begin{align}\label{equation:def_scale_dgff}
	\psi^N_v\coloneqq \int_{0}^{1} \sigma(s) \nabla \phi^N_v(s) \mathrm{d}s.
	\end{align}
	For $\delta>0$, let $V_N^\delta= [\delta N, (1-\delta)N)^2\cap \mathbb{Z}^2.$
 \cite[Lemma~3.3~(ii)]{paper2} shows that it is a centred Gaussian field with covariance given by
	\begin{align}
	\mathbb{E}\left[\psi^{N}_v\psi^N_w\right]= \log N \mathcal{I}_{\sigma^2}\left(\frac{\log N - \log_{+} \|v-w\|_2}{\log N}\right) +O(1),\quad \text{for } v,w\in V_N^\delta,
	\end{align}
	with $\log_{+} = \max \left\{0,\log(x)\right\}$.
\end{definition}
\begin{assumption}\label{assumption:1}
	In the rest of the paper, $\{\psi^N_v\}_{v\in V_N}$ is always a 2d scale-inhomogeneous DGFF on $V_N$. Moreover, we assume that $\mathcal{I}_{\sigma^2}(x)<x$, for $x\in (0,1),$ and that $\mathcal{I}_{\sigma^2}(1)=1$, with $s\mapsto \sigma(s)$ being differentiable at $0$ and $1$, such that $\sigma(0)<1$ and $\sigma(1)>1$.
\end{assumption}
Under \autoref{assumption:1} we proved in \cite{paper1,paper2}, building on work by  Arguin and Ouimet \cite{MR3541850}, the sub-leading order correction, tightness and convergence of the appropriately centred maximum. More explicitely, there exists a constant, $\beta=\beta(\sigma(1))$, which depends only on the final variance $\sigma(1)$, and a random variable, $Y=Y(\sigma(0))$, depending only on the initial variance $\sigma(0)$, such that, for any $z\in \mathbb{R}$,
\begin{align}\label{equation:convergence_max}
	\lim\limits_{N\rightarrow \infty}\mathbb{P}\left(\max_{v \in V_N} \psi^N_v\leq m_N-z\right)=\mathbb{E}\left[\exp\left[-\beta Y e^{-2z}\right]\right],
\end{align}
where $m_{ N}\coloneqq 2\log N - \frac{\log \log N}{4}.$
In particular, the limiting law solely depends on $\sigma(0)$ and $\sigma(1)$ and is therefore universal in the considered regime. Note that $m_N$ is also the maximum of $N^2$ i.i.d. $\mathcal{N}(0,\log N)$. Moreover, we proved in \cite[Theorem~2.2]{paper2} that under \autoref{assumption:1}, points whose height is close to the maximum are either $O(N)$ apart or within distance $O(1)$. In particular, there is a constant $c>0$, such that
\begin{align}\label{equation:1.7}
	\lim\limits_{r\rightarrow\infty}\lim\limits_{N\rightarrow \infty}\mathbb{P}\left(\exists u,v \in V_N\text{ with } r\leq \|u-v\|_2\leq\frac{N}{r}\text{ and }\psi^N_u,\psi^N_v\geq m_{ N}-c\log \log r \right)=0.
\end{align}
To state our results, we introduce some additional notation. Let $A\subset [0,1]^2$ and $B\subset \mathbb{R}$ be two Borel sets. For $v\in \mathbb{Z}^2$ and $r>0$, let its $r-$neighbourhood be $\Lambda_r(v)=\{w\in \mathbb{Z}^2: \|v-w\|_1\leq r \}$. Then, define
\begin{align}\label{equation:point_process_extremes}
	\eta_{N,r}(A\times B)\coloneqq \sum_{v\in V_N}\mathbbm{1}_{\psi^N_v=\max_{u\in \Lambda_r(v)}\psi^N_u}\,\mathbbm{1}_{x/N\in A}\mathbbm{1}_{\psi^N_v-m_N\in B}.
\end{align}
$\eta_{N,r}$ is a point measure encoding both position and relative height of extreme local maxima in $r-$neighbourhoods. To study distributional limits of these point measures, we equip the space of point measures on $[0,1]^2\times \mathbb{R}$ with the vague topology.
\begin{thm}\label{theorem:convergence_local_extremes}
	Let $\{\psi^N_v \}_{v\in V_N}$ be a scale-inhomogeneous DGFF satisfying \autoref{assumption:1}.
	Then, there is a random measure $Y(\mathrm{d}x)$ on $[0,1]^2$ which depends only on the initial variance $\sigma(0)$ and satisfies almost surely $Y([0,1]^2)<\infty$ and $Y(A)>0$, for any open and non-empty $A\subset [0,1]^2$. Moreover, there is a constant $\beta=\beta(\sigma(1))>0$, depending only on the final variance $\sigma(1)$, such that, for any sequence $r_N$ with $r_N\rightarrow \infty$ and $r_N/N\rightarrow 0$, as $N\rightarrow \infty$,
	\begin{align}
		\eta_{N,r_N} \overset{N\rightarrow \infty}{\rightarrow} PPP\left( Y(\mathrm{d}x)\otimes \beta e^{-2 h}\mathrm{d}h \right),
	\end{align}
	where convergence is in law with respect to the vague convergence of Radon measures on $[0,1]^2\times \mathbb{R}$.
\end{thm}
As the field at nearby vertices is strongly correlated, around each local maximum there will naturally be plenty of particles being close to it. Together with location and height of $r-$local maxima, we encode them in the point process
\begin{align}
\mu_{N,r}\coloneqq \sum_{v\in V_N}\mathbbm{1}_{\psi^N_v=\max_{u\in \Lambda_r(v)}\psi^N_u}\,\delta_{x/N}\otimes\delta_{\psi^N_v-m_N}\otimes\delta_{\{\psi^N_v-\psi^N_{v+w}:\, w\in \mathbb{Z}^2\}}.
\end{align}
These are Radon measures on $[0,1]^2\times \mathbb{R}\times \mathbb{R}^{\mathbb{Z}^2}$. We consider this space equipped with the topology of vague convergence. 
The following theorem shows convergence of $\mu_{N,r}$, the full extremal process.
\begin{thm}\label{thm:full_convergence}
	There is a probability measure $\nu$ on $[0,\infty)^{\mathbb{Z}^2}$ such that for each $r_N$ with $r_N\rightarrow \infty$ and $r_N/N\rightarrow 0$, as $N\rightarrow \infty$,
	\begin{align}
		\mu_{N,r_N}\rightarrow PPP\left( Y(\mathrm{d}x)\otimes \beta e^{-2 h}\mathrm{d}h\otimes \nu(\mathrm{d}\theta) \right).
	\end{align}
	The convergence is in law with respect to the vague convergence of Radon measures on $[0,1]^2\times \mathbb{R}\times \bar{\mathbb{R}}^{\mathbb{Z}^2}$. Moreover, $\nu$ is given by the weak limit,
	\begin{align}\label{equation:clusterlaw}
		\nu(\cdot)=\lim\limits_{r\rightarrow \infty}\mathbb{P}\left(\phi^{\mathbb{Z}^2\setminus\{0\}}+2\sigma(1)\mathfrak{a}\in \cdot|\phi^{\mathbb{Z}^2\setminus\{0\}}_w +2\sigma(1)\mathfrak{a}(w)\geq 0,\,\forall \|w\|_1\leq r \right),
	\end{align}
	with $\mathfrak{a}(w)=\lim\limits_{N\rightarrow \infty} G_{V_{2N}}\left[\left(N,N\right),\left(N,N\right)\right]-G_{V_{2N}}\left[\left(N,N\right),\left(N,N\right)+w\right]$ being the potential kernel. In addition, $\theta_0=0$ and $|\{w\in \mathbb{Z}^2:\,\theta_w \leq c\}|<\infty$, $\nu-$a.s. for each $c>0$.
\end{thm}
As a consequence of \autoref{thm:full_convergence}, we obtain convergence of the extremal process
\begin{align}
	\eta_N\coloneqq \sum_{v\in V_N} \delta_{v/N}\otimes \delta_{\psi^N_v-m_N}.
\end{align}
\begin{cor}
	Let $\{(x_i,h_i):\, i\in \mathbb{N} \}$ enumerate the points in a sample of $PPP\left( Y(\mathrm{d}x)\otimes \beta e^{-2h}\mathrm{d}h\right).$ Let $\{\theta^{(i)}_w:\, w\in \mathbb{Z}^2 \}$, $i\in \mathbb{N}$, be independent samples from the measure $\nu$, independent of $\{(x_i,h_i):\, i\in \mathbb{N} \}$. Then, as $N\rightarrow \infty$,
	\begin{align}\label{equation:1.113}
		\eta_N \rightarrow \sum_{i\in \mathbb{N}}\sum_{w\in \mathbb{Z}^2}\delta_{(x_i,h_i-\theta^{(i)}_w)}.
	\end{align}
	The convergence is in law with respect to the vague convergence of Radon measures on $[0,1]^2 \times \mathbb{R}$. Moreover, the measure on the right-hand side of \eqref{equation:1.113} is locally finite on $[0,1]^2\times \mathbb{R}$ a.s.
\end{cor}
	\subsection{Related work}
	Choosing $\sigma(x)\equiv 1$, for $x\in [0,1]$, in \eqref{equation:def_scale_dgff} gives the 2d DGFF.
Its maximum value was investigated by Bolthausen, Bramson, Daviaud, Deuschel, Ding, Giacomin and Zeitouni \cite{MR1880237,MR2243875,MR2772390,MR2846636,MR3101848,MR3262484,MR3433630}, which culminated in the proof of convergence of the maximum \cite{MR3433630}. Biskup and Louidor proved convergence of the extremal point process encoding local maxima and the field centred at those, to a cluster Cox process \cite{MR3509015,2016arXiv160600510B}. The random intensity measure is identified with the so-called Liouville quantum gravity measure \cite{MR4029145}. The cluster law of the 2d DGFF admits a closely related formulation to the one we obtain in \autoref{thm:full_convergence}, namely
\begin{align}\label{equation:1.13}
	\nu_{DGFF}=\lim\limits_{r\rightarrow \infty}\mathbb{P}\left(\phi^{\mathbb{Z}^2\setminus\{0\}}+2\mathfrak{a}\in \cdot|\phi^{\mathbb{Z}^2\setminus\{0\}}_w +2\mathfrak{a}(w)\geq 0,\,\forall \|w\|_1\leq r \right).
\end{align}
The slight, however important difference, is that the factor $\sigma(1)$ in \eqref{equation:clusterlaw} is equal to one. This causes the conditioning in \eqref{equation:1.13} to be asymptotically singular.
There is another possible regime in the scale-inhomogeneous DGFF, i.e. when $\mathcal{I}_{\sigma^2}(x)>x$, for some $x\in (0,1)$. When  $x\mapsto \mathcal{I}_{\sigma^2}(x)$ is piecewise linear, the leading and sub-leading order of the maximum, as well as exponential tails of the centred maximum, in particular tightness, are known \cite{MR3541850,paper1}. \par
Variable-speed branching Brownian motion (BBM), which first appeared in a paper by Derrida and Spohn \cite{MR971033}, is the natural analogue in the  context of BBM of the scale-inhomogeneous DGFF. It is a centred Gaussian process indexed by the leaves of the super-critical Galton-Watson tree, and covariance given by $tA(d(v,w)/t)$, where $d(v,w)$ is the time of the most recent common ancestor of two leaves $v$ and $w$. $A(x)\equiv 1$ corresponds to standard BBM. Its extremal process was investigated in \cite{MR3101852,MR3129797,2014arXiv1412.5975B,MR0494541,MR893913,MR2838339,MR2985174,MR3980921}. In \cite{MR3101852,MR3129797}, the cluster process was shown  to be BBM conditioned on the maximum being larger than $\sqrt{2}t$, or alternatively given as the limiting distribution of the neighbours of a local maximum. The extremal process of variable-speed BBM was investigated in \cite{MR3164771,MR3351476,MR3531703,MR2981635,1808.05445}. In the regime of weak correlations, i.e. when $A(x)<x$, for $x\in (0,1)$, $A^\prime (0)<1$ and $A^\prime(1)>1$, Bovier and Hartung \cite{MR3164771,MR3351476} proved convergence of the extremal process to a cluster Cox process. The cluster law can be described by the law of BBM in time $t$, conditioned on the maximum being larger than $\sqrt{2}A^\prime(1)t$, which is a perfect match to the one in the weakly correlated regime of the scale-inhomogeneous DGFF in \eqref{equation:clusterlaw}. In the regime when $A$ is strictly concave, Bovier and Kurkova \cite{MR2070335} showed that the first order of the maximum depends only on the concave hull of $A$. Moreover, Maillard and Zeitouni \cite{MR3531703} proved that the $2$nd order correction is proportional to $t^{1/3}$. \par
Note that there are other models such as the BRW \cite{MR3852245} or first passage percolation \cite{kistler2018oriented} where it was proven that the extremal process converges to a (cluster) Cox process.
	\subsection{Outline of Proof}
	We start to explain the proof of \autoref{theorem:convergence_local_extremes}. First, we deduce tightness of $\eta_{N,r}$ from \eqref{equation:convergence_max}, \eqref{equation:1.7} and a uniform exponential upper bound on extreme level sets, which is proven in \autoref{proposition:level_set_bound}. Then, we characterize possible limit laws as a Cox processes using a superposition principle as in \cite{MR3509015}. Finally, we need to show uniqueness of the random intensity measure. This follows from the convergence in distribution of multiple local maxima over disjoint subsets (see \autoref{theorem:conv_multiple_maxima}).
\par Next, we explain the proof of \autoref{thm:full_convergence}. By \eqref{equation:1.7}, we know that extreme local maxima have to be separated at distance $O(N)$ and, due to correlations, are surrounded by $O(1)$ neighbourhoods of high points. We need to show that the $O(1)$ neighbourhoods of extreme local maxima converge to independent samples of a cluster law. Using \eqref{equation:1.7} we know that also the $O(1)$ neigbourhoods must be at macroscopic distance, i.e. at distance of $O(N)$. To obtain independence of the clusters, we decompose the field into a sum of independent ``local fields'' that are zero outside the $O(1)$ neighbourhoods and a ``binding field'', which captures the contributions from outside the neighbourhoods. The requirement of being a cluster around a local maximum then translates into the local field being smaller than the value at its centre. We then show convergence of the laws of the local fields conditioned on a local maximum at their centre. In particular, we deduce that the clusters are i.i.d. samples of a common cluster law. Together with convergence of the extremal process of local maxima, \autoref{theorem:convergence_local_extremes}, this yields \autoref{thm:full_convergence}.
\par
\textit{Structure of the paper:}
In \autoref{section:pf1}, we prove \autoref{theorem:convergence_local_extremes}. The necessary ingredient, convergence of multiple local maxima over disjoint subsets, i.e. \autoref{theorem:conv_multiple_maxima}, is proved in \autoref{section:pf2}.
The proof of \autoref{thm:full_convergence} is provided in \autoref{section:pf3}. The appendix recalls Gaussian comparison tools.

	\section{Proof of \autoref{theorem:convergence_local_extremes}}\label{section:pf1}
	It turns out that we are able to follow and use large parts of the proof for the DGFF by Biskup and Louidor \cite{MR3509015}. As depicted in \cite{MR3509015,MR3618123}, the fact that the limiting point process takes the particular form of a generalized Poisson point process, is a consequence of a superposition property, which is due to its Gaussian nature along with certain properties of the field such as the separation of local maxima \cite{paper2} and tightness of extreme level sets. The main ingredient we need, in order to apply the machinery from \cite{MR3509015} to obtain the distributional invariance and thus Poisson limit laws, is tightness of the point processes, which is a consequence of the following proposition and previous results in \cite{paper2}.
For $y\in \mathbb{R}$, we denote by
\begin{align}
	\Gamma_N(y)=\left\{v\in V_N:\, \psi^N_v\geq m_N-y \right\},
\end{align}
the level set above $m_N-y$.
\begin{prop}\label{proposition:level_set_bound}
	There exists a constant $C>0,$ such that, for all $z>1$ and all $\kappa$,
	\begin{align}\label{equation:summable_levelset}
		\sup_{N\geq 1} \mathbb{P}\left(| \Gamma_N(y)| > e^{\kappa z}\right) \leq C e^{2y-\kappa z}.
	\end{align} 
	\begin{proof}
		By a first order Chebychev inequality and a standard Gaussian tail bound,
		\begin{align}
			\mathbb{P}\left(| \Gamma_N(y)| > e^{\kappa z}\right)\leq \tilde{C} \frac{\sqrt{\log N }}{m_N-\lambda} N^2 \exp\left[-\frac{(m_N-y)^2}{2\log N}\right] \leq C \exp\left[2y-\kappa z\right],
		\end{align}
		which shows \eqref{equation:summable_levelset}.
	\end{proof}
\end{prop}
\autoref{proposition:level_set_bound} together with \cite[Theorem~2.1]{paper2} implies tightness of $\{\eta_{N,r_N}\}_{N\in \mathbb{N}}$, as the right-hand side of \eqref{equation:summable_levelset} tends to zero as $N\rightarrow \infty$
	\subsection{Distributional Invariance}
	Let $(W_t)_{t\geq 0}$ be an independent standard Brownian motion started in $0$. Given a measurable function $f:[0,1]\times \mathbb{R}\rightarrow [0,\infty)$, let

\begin{align}
	f_t(x,h)=-\log \mathbb{E}^0\left[e^{-f(x,h+W_t-\frac{1}{2}t)} \right],\quad t\geq 0,
\end{align}
where $\mathbb{E}^0$ is the expectation with respect to the Brownian motion $(W_t)_{t\geq 0}$.
\begin{thm}(cp. \cite[Theorem~3.1]{MR3509015}\label{theorem:distributional_invariance})
	Let $\eta$ be any sub-sequential distributional limit of the processes $\{\eta_{N,r_N}\}_{N\geq 1}$, for some $r_N\rightarrow \infty$ with $r_N/N\rightarrow 0$. Then, for any continuous function $f:[0,1]^2\times \mathbb{R}\rightarrow [0,\infty)$ with compact support and all $t\geq 0$,
	\begin{align}
		\mathbb{E}\left[e^{-<\eta,f>}\right]=\mathbb{E}\left[e^{-<\eta, f_t>}\right].
	\end{align}
	\begin{proof}
		The proof of \autoref{theorem:distributional_invariance} is a rerun of the one in the case of the 2d DGFF \cite[Theorem~3.1]{MR3509015}. We therefore omit details here. It essentially uses convergence of the maximum obtained in \cite{paper2} together with expontential bounds on level sets, see \autoref{proposition:level_set_bound}.
	\end{proof}
\end{thm}
\begin{rem}
	As we think that the interpretation of the statement by Biskup and Louidor in \cite{MR3509015} is enlightening, we reproduce it here.
	Picking a sample, $\eta$, of the limit process, we know by tightness that $\eta(C)<\infty$ almost surely for any compact $C$. This allows us to write 
	\begin{align}
	\eta= \sum_{i \in \mathbb{N}}\delta_{(x_i,h_i)},
	\end{align}
	where $ \{(x_i,h_i)\in [0,1]\times \mathbb{R}\cup \{-\infty \}: i\in \mathbb{N}\}$ enumerate the points.
	Let $\{W^{(i)}_t:\, i\in \mathbb{N}\}$ be a collection of independent standard Brownian motions, independent of $\eta$, and set
	\begin{align}
	\eta_t \coloneqq \sum_{i \in \mathbb{N}} \delta_{(x_i,h_i+ W^{(i)}_t-\frac{1}{2}t)},\quad t\geq 0.
	\end{align}
	Using Fubini and dominated convergence, we have for all non-negative functions $f$,
	\begin{align}
	\mathbb{E}\left[e^{-<\eta, f_t>} \right]=\mathbb{E}\left[e^{-<\eta_t,f>} \right].
	\end{align}
	\autoref{theorem:distributional_invariance} then implies,
	\begin{align}
	\eta_t \overset{d}{=} \eta,\quad t\geq 0.
	\end{align}
\end{rem}
We borrow from \cite{MR3509015} a short heuristic argument why \autoref{theorem:distributional_invariance} should hold.
Let $\psi$ be a scale-inhomogeneous DGFF on $V_N$ satisfying \autoref{assumption:1} and let $\psi^\prime, \psi^{\prime \prime}$ be two independent copies of it. Fix some $t>0$.
Then,
\begin{align}
\psi\overset{d}{=} \sqrt{1-\frac{t}{\log N}} \psi^\prime +\sqrt{\frac{t}{\log N}}\psi^{\prime \prime}
=\psi^\prime -\frac{t}{2\log N} \psi^{\prime} +\sqrt{\frac{t}{\log N}} \psi^{\prime \prime} +o(1),
\end{align}
where we have used a Taylor expansion of the first square root, which has an error term $O(t^2/\log^2 N)$. Using the fact, that the first order of the maximum of the scale-inhomogeneous DGFF is $\log N$, we obtain an error $o(1)$. If we take $v\in V_N$ away from the boundary, where $\psi_v \geq m_N -y$ or $\psi^\prime_v\geq m_N -y$ and consider the $r-$neighbourhood $\Lambda_r(v)$, we first note that, for $w\in \Lambda_r(v)$, $\psi^{\prime \prime}_w -\psi_v^{\prime \prime}=O(1)$, and so by the prefactor, we may write,
\begin{align}
\psi_w \overset{d}{=} \psi^\prime_w -\frac{t}{2\log N} \psi^\prime_w + \sqrt{\frac{t}{\log N}} \psi^{\prime \prime}_v + o(1),\quad w\in \Lambda_r(v).
\end{align}
Similarly, we know that $\psi_w - m_N=O(1)$ and $\psi^\prime_w-m_N=O(1)$, for $w\in \Lambda_r(v)$, and thus, we may replace $\frac{t}{2\log N} \psi^\prime_w$ by $\frac{t}{2\log N}(m_N+O(1))=t+o(1)$, to obtain
\begin{align}
\psi_w \overset{d}{=} \psi^{\prime}_w - t + \sqrt{\frac{t}{\log N}} \psi^{\prime \prime}_v +o(1),\quad w\in \Lambda_r(v).
\end{align}
Finally, we see that $\sqrt{\frac{t}{\log N}}\psi^{\prime \prime}$ is asymptotically distributed as $W_t$, where $(W_t)_{t\geq 0}$ is a Brownian motion. Further, we know from \cite[Theorem~2.2]{paper2}, that local extremes are at distance of order $N$ and so the field $\psi^{\prime \prime}$ in two such neighbourhoods has correlation of order $O(1)$. The normalizing factor $\sqrt{\frac{t}{\log N}}$ then implies that two such neighbourhoods are asymptotically independent. Thus, for $N$ large, we have a one-to-one correspondence between local maxima of $\psi$ and local maxima of $\psi^\prime$ by a shift in their height through independent Brownian motions with drift $-1$.
	\subsection{Poisson limit law}
	Just as in \cite{MR3509015}, distributional invariance, \autoref{theorem:distributional_invariance}, allows to extract a Poisson limit law for every such subsequence, i.e. for any sub-sequential limit of the extremal process. In our setting, we can directly apply \cite[Theorem~3.2]{MR3509015}.
\begin{thm}\cite[Theorem~3.2]{MR3509015}\label{theorem:poisson_limit_law}
	Suppose that $\eta$ is a sub-sequential limit of the process $\eta_{N,r_N}$, that is a point process on $[0,1]^2\times \mathbb{R}$ such that, for some $t>0$, and all continuous functions $f: [0,1]^2 \times \mathbb{R}\rightarrow [0,\infty)$ with compact support, it holds, as in \autoref{theorem:distributional_invariance},
	\begin{align}
		\mathbb{E}\left[e^{-<\eta,f>}\right]=\mathbb{E}\left[e^{-<\eta,f_t>}\right].
	\end{align}
	Moreover, assume that almost surely $\eta([0,1]^2\times [0,\infty))<\infty$ and $\eta([0,1]^2\times \mathbb{R})>0$. Then, there is a random Borel measure $Y$ on $[0,1]^2$, satisfying $Y([0,1]^2)\in (0,\infty)$ almost surely, such that
	\begin{align}
		\eta \overset{d}{=} PPP\left(Y(\mathrm{d}x)\otimes \beta e^{-2 h}\mathrm{d}h\right).
	\end{align}
\end{thm}
	\subsection{Uniqueness}
	In this section, we show uniqueness of the extremal process of local extremes, i.e. of the limit $\lim\limits_{N\rightarrow \infty}\eta_{N,r_N}$. In light of \autoref{theorem:poisson_limit_law}, we do this by showing uniqueness of the random measure $Y(\mathrm{d}x)$. The proof is a generalization of the proof of uniqueness of the random variable $Y$ in \cite[Theorem~2.1]{paper2}. We show that the joint law of local maxima converges in law and that this law can be written as a Laplace transform of the random measure $Y(\mathrm{d}x)$, which then implies uniqueness of $Y(\mathrm{d}x)$. For a set $A\subset [0,1]$, we write $\psi^{*}_{N,A}=\max\left\{\psi^N_v:\,v\in V_N,\, v/N\in A \right\}$.
\begin{thm}\label{theorem:conv_multiple_maxima}
	Let $(A_1, \dotsc,A_p)$ be a collection of disjoint non-empty open subsets of $[0,1]^2$. Then the law of $\left(\max\{\psi^N_v:\, v\in V_N,\,v/N \in A_l \} -m_N \right)_{l=1}^{p}$ converges weakly as $N\rightarrow \infty$. More precisely, there are random variables $Y_{A_1},\dotsc, Y_{A_p} $ depending only on the initial variance $\sigma(0)$, satisfying $Y_{A_i}>0$ almost surely, for $1\leq i \leq p$, and there is a constant $\beta>0$, depending only on the final variance $\sigma(1)$, such that
	\begin{align}
		\lim\limits_{N\rightarrow \infty} \mathbb{P}\left(\psi^{*}_{N,A_l}-m_N\leq x_l:\, l=1,\dotsc,p \right)= \mathbb{E}\left[\exp\left(-\beta  \sum_{l=1}^pe^{-2 x_l} Y_{A_l} \right) \right].
	\end{align}
\end{thm}
The constant $\beta$ in \autoref{theorem:conv_multiple_maxima} is identical to the one appearing in \eqref{equation:convergence_max}. Next, we prove \autoref{theorem:convergence_local_extremes}. The proof of \autoref{theorem:conv_multiple_maxima} is given in \autoref{section:pf2}. 
	\begin{proof}[Proof of \autoref{theorem:convergence_local_extremes} using \autoref{theorem:conv_multiple_maxima}]
	Let $r_N\rightarrow \infty$ with $r_N/N\rightarrow 0$ be now a fixed sequence. Denote by $\eta$ a corresponding sub-sequential limit of the extremal process $\{\eta_{N,r_N}\}_{N\geq 1}$. By \autoref{theorem:poisson_limit_law}, there is a corresponding random measure $\tilde{Y}(\mathrm{d}x)$ such that $\eta \overset{d}{=} PPP\left(\tilde{Y}(\mathrm{d}x)\otimes \beta e^{-2 h}\mathrm{d}h\right).$
	Note that, as a trivial consequence of \autoref{theorem:conv_multiple_maxima}, for any open and non-empty $A\subset [0,1]^2$, $\psi^{*}_{N,A}-m_N$ is a tight sequence. Fix an arbitrary collection, $(A_1,\dotsc, A_p)$, of disjoint, open and non-empty subsets of $[0,1]^2$, with $\tilde{Y}(\partial A_l)=0$, for any $l\in \{1,\dotsc,p\}$.
	By \autoref{theorem:conv_multiple_maxima}, there is a dense subset $R\subset \mathbb{R}$ such that, for any $x_1,\dotsc,x_p\in R$,
	\begin{align}\label{equation:lapalce_trafo_multiple_maxima_disjoint_subsets}
		\mathbb{E}\left[\exp\left(-\beta \sum_{l=1}^{p}e^{-2x_l}\tilde{Y}(A_l) \right)\right]= \lim\limits_{N\rightarrow \infty} \mathbb{P}\left(\psi^{*}_{N,A_l}-m_N\leq x_l:\, l=1,\dotsc,p \right).
	\end{align}
	Again by \autoref{theorem:conv_multiple_maxima}, the right-hand side of \eqref{equation:lapalce_trafo_multiple_maxima_disjoint_subsets} is the same for all subsequences. Using continuity in $x$ of the left hand side, we can deduce from convergence on the dense subset $R$, convergence on $\mathbb{R}$. Along with a standard approximation argument of continuous functions on $[0,1]^2$ via non-negative simple functions, this implies uniqueness of the Laplace transform of the random measure $\tilde{Y}(\mathrm{d}x)$ on the disjoint collection $(A_1,\dotsc,A_p)$, regardless of the subsequence considered. As $p\in \mathbb{N}$ and $A_1,\dotsc, A_p$ are arbitrary, it follows that $\tilde{Y}(\mathrm{d}x)$ is the same for all sub-sequences.
	Therefore, we obtain a random Borel measure $Y(\mathrm{d}x)$ whose masses of any countable collection of open sets $A_1,\dotsc, A_p$ are given by $Y_{A_1},\dotsc, Y_{A_p}$ from \autoref{theorem:conv_multiple_maxima}, depending only on $\sigma(0)$. We conclude, that the law of the measure $Y(\mathrm{d}x)$ also depends only on initial variance, $\sigma(0)$.
	Further, note that by \autoref{proposition:level_set_bound},
	\begin{align}\label{equation:2.18}
		\mathbb{P}\left(\eta([0,1]^2\times [-y,\infty])>e^{k y}\right)\leq Ce^{-y(\kappa-2)}.
	\end{align}
	 In combination with \autoref{theorem:poisson_limit_law}, \eqref{equation:2.18} implies that the total mass of $Y$ is almost surely finite. Moreover, \autoref{theorem:conv_multiple_maxima} implies that, for any non-empty and open $A\subset [0,1]^2$, we have almost surely $Y(A)>0$.
\end{proof}
	%\section{Proof of \autoref{theorem:size_of_level_sets}}\label{section:proof_tightness_level_sets}
	%\input{proof_dist_inv.tex}
	%\subsection{Proof of \autoref{theorem:distributional_invariance}}
	%\input{proof_dist_invariance.tex}
	%\section{Poisson limit law}
	%\input{proof_poisson_limit.tex}
		%\section{Proof of \autoref{theorem:distributional_invariance}}
	%\input{dist_invariance.tex}
	\section{Proof of \autoref{thm:full_convergence}}\label{section:pf3}
	In the following, we assume that $V_N$ is centred at the origin. Let $\mu$ be a Radon measure on $[0,1]^2\times \mathbb{R} \times \mathbb{R}^{\mathbb{Z}^2}$ and $f: [0,1]^2\times \mathbb{R} \times \mathbb{R}^{\mathbb{Z}^2}\rightarrow [0,\infty)$ be a measurable function with compact support. We write
\begin{align}
	\langle \mu, f\rangle \coloneqq \int \mu(\mathrm{d}x\mathrm{d}h\mathrm{d}\theta)f(x,h,\theta).
\end{align}
Further, let
\begin{align}
	\Theta_{N,r}\coloneqq \{v\in V_N:\, \psi^N_v=\max_{u\in \Lambda_r(v)} \psi^N_u \}
\end{align}
be the set of $r-$local maxima.
\begin{lemma}\label{lemma:1}
	For any $r_N\rightarrow \infty$ with $r_N/N\rightarrow 0$ and any continuous function $f: [0,1]^2\times \mathbb{R}\times \mathbb{R}^{\mathbb{Z}^2}$ with compact support,
	\begin{align}\label{equation:lemma1}
		\lim\limits_{r\rightarrow \infty}\limsup\limits_{N\rightarrow \infty}\max_{M:r\leq M\leq N/r}\left|\mathbb{E}\left[ e^{-\langle \mu_{N,r_N}}\rangle\right]-\mathbb{E}\left[ e^{-\langle \mu_{N,M},f\rangle}\right] \right|=0.
	\end{align}
	\begin{proof}
		Let $\lambda>0$ be such that $f(x,h,\theta)=0$, for $h\geq \lambda$. If $\langle \mu_{N,r_N},f\rangle\neq \langle \mu_{N,M},f\rangle$, for some $M$ with $r\leq M\leq N/r$, then $\Theta_{N,r_N}\triangle \Theta_{N,
		M}\cap \Gamma_N(\lambda)\neq \emptyset.$ Thus, there are $u,v\in \Gamma_N(\lambda)$ such that $\min(M,r_N)\leq \|u-v\|_2\leq \max(M,r_N)$.
		For $N$ being so large that $ r_N>r$ and $r_N\leq N/r$, this implies
		\begin{align}
			\max_{M:r\leq M\leq N/r}\left| \mathbb{E}\left[e^{-\langle \mu_{N,r_N},f\rangle}\right]-\mathbb{E}\left[e^{-\langle \mu_{N,M},f\rangle}\right]\right|\leq \mathbb{P}\left(\exists u,v\in \Gamma_N(\lambda):\, r\leq \|u-v\|_2\leq N/r \right),
		\end{align}
		which by \cite[Theorem~2.2]{paper2} tends to zero. This shows \eqref{equation:lemma1}.
	\end{proof}
\end{lemma}
We set $M\coloneqq \min \{k: 2^k> r\}$. In light of \autoref{lemma:1}, we work with $\mu_{N,M}$ instead of $\mu_{N,r_N}$.
Suppose that the local maximum is taken at $v\in V_N$. We decompose into two fields. The idea is, for fixed $v\in V_N$, to use the Gibbs-Markov property of the underlying DGFF to write the field into independent components. One that captures the field inside $\Lambda_{M}(v)$ and another that captures the field outside, i.e. in $\Lambda_{M}^c(v)$. $v\in V_N$ later plays the role of a local maximum. Thus, we write
\begin{align}\label{equation:3..5}
\psi^N_w= \Phi^{M,v}_w+\tilde{\psi}^{\Lambda_{M}(v)}_w,\quad \text{ for } w\in \Lambda_{M}(v),
\end{align}
where
\begin{align}\label{equation:3.5}
	\Phi^{M,v}_w\coloneqq \int_{0}^{1-\frac{\log M+log_{+}\|v-w\|_2}{\log N}}\sigma(s)\nabla\phi^N_w(s)\mathrm{d}s +\int_{1-\frac{\log M+log_{+}\|v-w\|_2}{\log N}}^{1}\sigma(s) \nabla \mathbb{E}\left[\phi^N_w| \sigma\left(\phi^N_y:y\in \partial [w]_s\cap \Lambda_{M}^c(v) \right)\right].
\end{align}
and where
\begin{align}\label{equation:3.8}
\tilde{\psi}^{\Lambda_M(v)}_w= \int_{1-\frac{\log M+\log_+ \|v-w\|_2}{\log N}}^{1}\sigma(s)\phi^{\Lambda_M(v)}_w(s)\mathrm{d}s.
\end{align}
The field in \eqref{equation:3.5} encodes the increments when conditioning outside the local maximum $v\in V_N$ and its $M-$neighbourhood, $\Lambda_{M}(v)$. The field in \eqref{equation:3.8} encodes the remaining increments within $\Lambda_M(v)$.
\begin{comment}
	If we further want to condition the field in \eqref{equation:3.8} on the value in its centre, $\tilde{\psi}^{\Lambda_M(v)}_v$, we can decompose it further into
	\begin{align}\label{equation:3.7}
	\tilde{\psi}^{\Lambda_M(v) \setminus\{v\}}_w=\int_{1-\frac{\log M+\log_+ \|v-w\|_2}{\log N}}^{1}\sigma(s)\nabla \phi^{\Lambda_M(v)\setminus\{v\}}_w(s)\mathrm{d}s,
	\end{align} 
	which is zero at $v$ and into
	\begin{align}\label{equation:3.6}
	\tilde{\psi}^{\Lambda_M(v)}_v\int_{1-\frac{\log M+log_+ \|v-w\|_2}{\log N}}^{1}\sigma(s)\nabla \left(\sum_{z \in \partial [w]_s}\mathbb{P}_w\left(
	S_{\tau_{\partial [w]_s}}=z\right)g_{M}(z-v)\right)\mathrm{d}s,
	\end{align}
	where $g_M(z-v)$ is the unique discrete harmonic being $1$ at the origin and $0$ on $\Lambda^c_M(0)$.
	The field in \eqref{equation:3.6} comes from conditioning on the local maximum value $\tilde{\psi}^{\Lambda_{M}(v)}_v$ and is the scale-inhomogeneous discrete harmonic extension of the centre value into the box $\Lambda_{M}(v).$
\end{comment} 
The following lemma points out the key idea behind the definitions in \eqref{equation:3.5} and \eqref{equation:3.8}.
\begin{lemma}\label{lemma:5.2}
	Suppose $v\in V_N$ such that $\Lambda_M(v)\subset V_N$ and let $M=2^k$. Consider the sigma-algebra
	\begin{align}
		\mathcal{F}_{M,v}\coloneqq \sigma\left(\phi^N_w:\, w\in \{v\}\cup \Lambda_M(v)^c\right).
	\end{align}
	Then, for Lebesgue almost every $t\in \mathbb{R}$,
	\begin{align}
		\mathbb{P}&\left(\psi^N_{v+\cdot}-\Phi^{M,v}_{v+\cdot}\in \cdot | \mathcal{F}_{M,v}\right)
		%&=\mathbb{P}\left(\tilde{\psi}^{\Lambda_M(v)\setminus\{v\}}_\cdot + \int_{1-\frac{\log M}{\log %N}}^{1}\sigma(s)\nabla \left(\sum_{z \in \partial [w]_s}\mathbb{P}_w\left(
		%S_{\tau_{\partial [w]_s}}=z\right)g_{M}(z-v)\right)\mathrm{d}s \in \cdot \right)\nonumber\\
		=\mathbb{P}\left(\tilde{\psi}^{\Lambda_M(v)}_{v+\cdot} \in \cdot | \tilde{\psi}^{\Lambda_M(v)}_v= t-\Phi^{M,v}_v\right),\qquad \text{on } \{\psi^N_v=t\}.
	\end{align}
	\begin{proof}
		It is an immediate consequence using \eqref{equation:3..5}.
	\end{proof}
\end{lemma}
The following proposition is used to localize the initial increments, $\Phi^{M,v}_v$, of a local maximum at $v\in V_N$.
\begin{prop}\label{prop:localization}
	Let $t\in \mathbb{R}$. There is $r_0\in \mathbb{N}$ such that, for any $\delta\in (0,1)$, $r\geq r_0$, $N\in \mathbb{N}$, sufficiently large, $M\in (r,N/r)$ and $\gamma\in (0,1/2)$, there is a constant $C_\delta>0$, depending only on $\delta$,
	\begin{align}
		\mathbb{P}\left(\exists v \in V_N:\, \psi^N_v\geq m_N- t, \Phi^{M,v}_v-2\log N \mathcal{I}_{\sigma^2}\left(1-\frac{\log M}{\log N}\right) \notin [-\log^\gamma (M),\log^\gamma (M)]\right)\nonumber\\
		\leq C_\delta e^{2s} \sum_{k=\lfloor \log M\rfloor}^{\infty} k^{\frac{1}{2}-\gamma}\exp\left[-k^{\frac{2\gamma-1}{2}}\right].
	\end{align}
	\begin{proof}
		As in \eqref{equation:3..5},
		\begin{align}\label{equation:3.18}
			\psi^N_v= \Phi^{M,v}_v+	\tilde{\psi}^{\Lambda_M(v)}_v,
		\end{align}
		where the fields on the right hand side are independent.
		Using \cite[Lemma~3.1~(i)]{paper2} for the first and the last field in \eqref{equation:3.18}, as well as by Green function asymptotics, see e.g. \cite[(3.47), (B.5)]{2016arXiv160600510B}, we deduce that, for any $\delta>0$, there is a constant $c_\delta>0$, such that
		\begin{align}
			\sup\limits_{v\in V^\delta_{N}}\mathrm{Var}\left[\Phi^{M,v}_v\right] \leq 2\log N \mathcal{I}_{\sigma^2}\left(1-\frac{\log M}{\log N}\right) +c_\delta.
		\end{align}
		Moreover, $\{\Phi^{M,v}_v\}_{v\in V_N}$ is a centred Gaussian field. Thus, we can rerun the proof of \cite[Proposition~4.2]{paper2}, where the constant on the right of \cite[(4.13)]{paper2} may now depend on $\delta$.
		This concludes the proof of \autoref{prop:localization}.
	\end{proof}
\end{prop}
The following lemma allows us to reduce the local field defined in \eqref{equation:3.8} to a usual DGFF with a constant parameter.
\begin{lemma}\label{lemma:3.4}
	Let $v\in V_N^\delta$ and let $\{\tilde{\psi}^{\Lambda_{M}(v)}_w:\, w\in \Lambda_{M}(v) \}$ be the centred Gaussian field defined in \eqref{equation:3.8}. Then,
	\begin{align}\label{equation:3.13}
		\lim\limits_{M\rightarrow \infty}\tilde{\psi}^{\Lambda_{M}(v)}-\sigma(1)\phi^{\Lambda_{M}(v)}=0 \quad a.s.
	\end{align}
	\begin{proof}
			Note that for some $\epsilon>0$, by an Taylor expansion at $s=1$, we have $\sigma(s)=\sigma(1)-\sigma^\prime (1)(1-s)+o(\sigma^\prime(1)(1-s))$, for $s\in (1-\epsilon,1]$. In particular, for any $v\in V_N$ and $w\in \Lambda_{M}(v)$,
		\begin{align}\label{equation:4.15}
		\tilde{\psi}^{\Lambda_{M}(v)}_{w}-\sigma(1)\phi^{\Lambda_M(v)}_w=\int_{1-\frac{\log M+\log_+ \|v-w\|_2}{\log N}}^{1} \sigma^\prime(1)(1-s)\nabla \phi^{\Lambda_{M}(v)}_w(s) \mathrm{d}s+o(1),
		\end{align}
		which is a centred Gaussian and where the error term vanishes, as $N\rightarrow \infty$. By Cauchy-Schwarz and asymptotics of the potential kernel, e.g. \cite[(2.7), (B.6)]{2016arXiv160600510B}, the covariances of the field on the right-hand side of \eqref{equation:4.15} is bounded by a uniform constant times $\log^2 M/ \log^{3/2} N$, which tends to zero uniformly, as $N\rightarrow \infty$. This shows \eqref{equation:3.13}
	\end{proof}
\end{lemma}
\begin{rem}\label{rem:1}
	With regard to \autoref{prop:localization}, the cluster law around around a local maximum $v\in V_N^\delta$ can be written in the form $\mathbb{P}\left(\tilde{\psi}^{\Lambda_{M}(v)}\in \cdot | \tilde{\psi}^{\Lambda_{M}(v)}_v=2\log N\mathcal{I}_{\sigma^2}\left(1-\frac{\log M}{\log N},1\right)+t, \tilde{\psi}^{\Lambda_M(v)}_w\leq \tilde{\psi}^{\Lambda_{M}(v)}_v\right)$. \autoref{lemma:3.4} shows that this has the same weak limit, as $M\rightarrow \infty$ after $N\rightarrow \infty$, as
	\begin{align}\label{equation:4.18}
	\nu^{(M,t)}(\cdot)\coloneqq\mathbb{P}\left(\sigma(1)\left(\phi^{\Lambda_M(0)}_0-\phi^{\Lambda_M(0)}\right)\in \cdot | \sigma(1)\phi^{\Lambda_M(0)}_0=2\sigma^2(1)\log M+t, \sigma(1)\phi^{\Lambda_M(0)}\leq \sigma(1)\phi^{\Lambda_M(0)}_0 \right).
	\end{align}
\end{rem}
In the following lemma we show that the the cluster limit of the law $\nu^{(M,t)}$ exists in a suitable sense.
\begin{lemma}\label{lemma:5.5}
	Fix $r,j\geq 1$ and let $c_1\in (0,\infty)$. For $M=\min\{k:\, 2^k >r \}$, uniformly in $f\in C_b\left(\mathbb{R}^{\Lambda_j}\right)$ and $t=o(\log M)$,
	\begin{align}
		\lim\limits_{M\rightarrow \infty}\mathbb{E}_{\nu^{(M,t)}}\left[f\right]= \mathbb{E}_\nu\left[f\right],
	\end{align}
	where $\nu(\cdot)\coloneqq \lim\limits_{r\rightarrow \infty}\nu_r(\cdot),$
	\begin{align}\label{equation:3.17}
		\nu_r(\cdot)\coloneqq \mathbb{P}\left(\phi^{\mathbb{Z}^2\setminus\{0\}}+2\sigma(1)\mathfrak{a}\in \cdot | \phi^{\mathbb{Z}^2\setminus\{0\}}_v+2\sigma(1)\mathfrak{a}(v)\geq 0:\, \|v\|_1 \leq r \right)
	\end{align}
	and $\mathfrak{a}$ being the potential kernel.
	\begin{proof}
		Convergence of the finite dimensional distributions of the measures $\nu_r(\cdot)$ is a simple consequence of the DGFF satisfying the strong FKG-inequality, which implies that $r\mapsto \nu_r$ is stochastically increasing. Thus, $\lim\limits_{r\rightarrow \infty}\nu_r(A)$ exists for any event $A$, depending on only a finite number of coordinates.
		Next, we prove that $\{\nu_r\}_r$ is tight, which then implies that $\nu$ is a distribution on $\mathbb{R}^{\mathbb{Z}^2}$. By a union and a Gaussian tail bound, for any $r\geq k_0>0$, there are constants $C,\tilde{C}>0$ such that
		\begin{align}
		\mathbb{P}\left(\exists v, k_0\leq \|v\|_1\leq r: \,\phi^{\mathbb{Z}^2\setminus\{0\}}_v> 2\sigma(1) \log \|v\| \right)\leq \sum_{k=k_0}^{r} 4k \mathbb{P}\left(\sup_{\|v\|_1=k} \phi^{\mathbb{Z}^2\setminus\{0\}}_v>2 \sigma(1)\log k+\frac{1}{2}\log(2) \right)\quad\nonumber\\ \leq C \sum_{k=k_0}^{r} \frac{4k}{\sqrt{\log k}}\exp\left[-\sigma^2(1)\log k+c_0\right]\leq \tilde{C} \sum_{k=k_0}^{\infty} \frac{1}{\sqrt{\log k}}\exp\left[-[\sigma^2(1)-1]\log k\right].
		\end{align}
		As the sum converges and vanishes, as $k_0\rightarrow \infty$, we deduce tightness of $(\nu_r)_{r\in \mathbb{N}}$ and so $\nu(\mathbb{R}^{\mathbb{Z}^2})=1.$
		In the last step, we show that it takes the particular form as in \eqref{equation:3.17}.
		We have that $\phi^{\Lambda_M(0)}$ conditioned on $\phi^{\Lambda_M(0)}_0=2\sigma(1)\log M$ shifts the mean of $\phi^{
			\Lambda_M(0)}_{0}-\phi^{\Lambda_M(0)}$ by a quantity with asymptotic
		\begin{align}
			(2\sigma(1)\log M +t)(1-g_M(v))\rightarrow 2\sigma(1)\mathfrak{a}(v),
		\end{align}
		as $M\rightarrow \infty$, and where $g_M(x)$ is discrete harmonic with $g_M(0)=1$ and $g_M(x)=0,$ for $x\notin \Lambda_{M}(0)$. In particular, the law of $v\mapsto \phi^{\Lambda_M(0)}_0-\phi^{\Lambda_M(0)}_v$ conditioned on $\phi^{\Lambda_{M}(0)}_0=2\sigma(1)\log M$ converges in the sense of finite dimensional distributions to
		\begin{align}
			\phi^{\mathbb{Z}^2\setminus \{0\}}_v+2\sigma(1)\mathfrak{a}(v),
		\end{align}
		where $\{\phi^{\mathbb{Z}^2\setminus \{0\}}_v\}_{v\in \mathbb{Z}^2\setminus \{0\}}$ is the pinned DGFF, which is a centred Gaussian field with covariances as in \cite[(2.7)]{2016arXiv160600510B}.
		This concludes the proof of \autoref{lemma:5.5}. 
	\end{proof}
\end{lemma}
Having weak convergence of the auxiliary cluster law, $\nu_r$, we are now in a position to prove convergence of the full extremal process.
\begin{proof}[Proof of \autoref{thm:full_convergence}]
	First note that by \autoref{lemma:1} we can work with $M$ instead of $r_N$. Let $f:[0,1]^2\times \mathbb{R}\times\mathbb{R}^{\mathbb{Z}^2}\mapsto [0,\infty)$ be a continuous function with compact support. In addition, assume that, for any $x\in [0,1]^2$ and $t\in \mathbb{R}$, $f(x,t,\phi)$ depends only on $\{\phi_y:\, y\in \Lambda_{M}(x)\}$. Let $V_N=\cup_{i=1}^{(N/M)^2}V_{M,i}$ be a decomposition of $V_N$ into disjoint shifts of $V_M$. Moreover, let $\delta\in (0,1)$ and set
	\begin{align}
		\mu_{N,M,\delta}\coloneqq\sum_{v\in \cup_{i=1}^{(N/M)^2} V_{M,i}^\delta}\mathbbm{1}_{v \in \Theta_{N,M}} \delta_{v/N}\otimes\delta_{\psi^N_v-m_N}\otimes\delta_{\{\psi^N_v-\psi^N_{v+w}:\, w\in \mathbb{Z}^2\}}.
	\end{align}
	By \autoref{prop:localization}, \cite[Proposition~5.1]{paper2} and \cite[Theorem~2.2]{paper2}, it suffices to compute
	\begin{align}\label{equation:4.23}
	\lim\limits_{\delta \rightarrow 0}\lim\limits_{M\rightarrow \infty}\lim\limits_{N\rightarrow \infty}\mathbb{E}\left[e^{-\langle \mu_{N,M,\delta},f\rangle}\mathbbm{1}_{N \|v-w\|_2>4 M:v,w \in \Theta_{N,M}}\mathbbm{1}_{\{\Phi^{N,v}_v-2\log N \mathcal{I}_{\sigma^2}\left(1-\frac{\log M}{\log N}\right)\in [-\log^\gamma (M),\log^\gamma (M)]:\, v \in \Theta_{N,M}\}}\right].
	\end{align}
	Set
	\begin{align}
	&	f_{N,M}(v/N,t)\coloneqq\nonumber\\& -\log \mathbb{E}\left[\exp\left[-f\left(x,t,\left(\psi^N_{v}-\Phi^{M,v}_v-\psi^N_{v+w}+\Phi^{M,v}_{v+w}:\,w\in \mathbb{Z}^2\right)\right)\right]| \psi^N_{v}=m_N+t, v\in \Theta_{N,M}\right].
	\end{align}
	Conditioning on position, $x_i N$, and height, $m_N+t_i$, of local maxima in $\cup_{i=1}^{(N/M)^2}V_{M,i}^{\delta}$ and on the sigma-algebra $\sigma\left(\phi^N_w:\,w\in \cup \partial\Lambda_M(x_i N)\right)$, using \autoref{lemma:5.2} and the Taylor approximation for the cluster process as in \autoref{rem:1}, we can rewrite \eqref{equation:4.23} as
	\begin{align}\label{equation:4.26}
		\mathbb{E}\left[\prod_{i=1}^{(N/M)^2}e^{-f_{N,M}(x_i,t_i)}\mathbbm{1}_{N \|x_j-x_k\|_2>4 M:x_j N, x_k N \in \Theta_{N,M} }\mathbbm{1}_{\{\Phi^{N,v}_v-2\log N \mathcal{I}_{\sigma^2}\left(1-\frac{\log M}{\log N}\in [-\log^\gamma (M),\log^\gamma (M)]\right):\, v \in \Theta_{N,M}\}}\right].
	\end{align}
	On $\{\Phi^{N,v}_v-2\log N \mathcal{I}_{\sigma^2}\left(1-\frac{\log M}{\log N}\right)\in [-\log^\gamma (M),\log^\gamma (M)]:\, v\in \Theta_{N,M}\}$, \autoref{lemma:5.2}, \autoref{lemma:3.4}, \autoref{rem:1} and \autoref{lemma:5.5} imply
	\begin{align}\label{equation:4.27}
		\lim\limits_{M\rightarrow \infty}\lim\limits_{N\rightarrow \infty} f_{N,M}(x,t)=f_{\nu}(x,t)\coloneqq - \log \mathbb{E}_{\nu}\left[e^{-f(x,t,\phi)}\right].
	\end{align}
	In particular, the convergence in \eqref{equation:4.27} is uniform in $x\in \cup_{i=1}^{(N/M)^2}V_{M,i}^{\delta}$ and $t\in \mathbb{R}$.
	Using \eqref{equation:4.26} and \autoref{prop:localization}, we can rewrite \eqref{equation:4.23} as
	\begin{align}\label{equation:4.28}
		\mathbb{E}\left[e^{-\langle \eta_{N,M},f_v\rangle}\right]+o(1).
	\end{align}
	Applying \autoref{theorem:convergence_local_extremes} to \eqref{equation:4.28}, we obtain
	\begin{align}\label{equation:5.27}
		\lim\limits_{M\rightarrow \infty}\lim\limits_{N\rightarrow \infty}\mathbb{E}\left[e^{-\langle \mu_{N,M},f_{\nu}\rangle}\right]&=\mathbb{E}\left[\exp\left[-\int_{[0,1]^2\times \mathbb{R}} Y(\mathrm{d}x)\otimes \beta e^{-2h}\mathrm{d}h\left(1-e^{-f_{\nu}(x,h)}\right)\right]\right]\nonumber\\
		&=\mathbb{E}\left[\exp\left[-\int_{[0,1]^2\times \mathbb{R}\times \mathbb{R}^{\mathbb{Z}^2}} Y(\mathrm{d}x)\otimes \beta e^{-2h}\mathrm{d}h\otimes \nu(\mathrm{d}\phi)\left(1-e^{-f(x,h,\phi)}\right)\right]\right].
	\end{align}
	Noting that the last line in \eqref{equation:5.27} is the Laplace transform of a Poisson point process with intensity $\beta Y(\mathrm{d}x)\otimes e^{-2h}\mathrm{d}h\otimes \nu(\mathrm{d}\phi)$, concludes the proof.
\end{proof}
	\section{Proof of \autoref{theorem:conv_multiple_maxima}}\label{section:pf2}
	First, we recall the $3-$field approximation used in \cite{paper2} to prove convergence in law of the centred maximum.

	\subsection{$3-$field approximation}
	We first decompose the underlying grid $V_N$. Assume $N=2^n$ to be much larger than any other forthcoming integers. Next, pick two large integers $L=2^l$ and $K=2^k$. Partition $V_N$ in a disjoint union of $(KL)^2$ boxes, $\mathcal{B}_{N/KL}=\{B_{N/KL,i}:\, i=1,\dotsc, (KL)^2 \}$, each of side length $N/KL$. Let $v_{N/KL,i}\in V_N$ be the left bottom corner of box $B_{N/KL,i}$ and write $w_{i}=\frac{v_{N/KL,i}}{N/KL}$. We consider $\{w_i\}_{i=1,\dotsc, (KL)^2}$ as the vertices of a box $V_{KL}$. Analogously, let $K^\prime=2^{k^\prime}$ and $L^\prime=2^{l^\prime}$ be two integers, such that $K^\prime L^\prime$ divides $N$. Let $\mathcal{B}_{K^\prime L^\prime}=\{B_{K^\prime L^\prime,i}:\, i=1,\dotsc,[N/(K^\prime L^\prime )]^2 \}$ be a disjoint partitioning of $V_N$ with boxes $B_{K^\prime L^\prime, j}$, each of side length $K^\prime L^\prime$. The left bottom corner of a box $B_{K^\prime L^\prime,i}$ we call $v_{K^\prime L^\prime,i}$.\\
We take limits in the order $N,L,K,L^\prime$ and then $K^\prime$, for which we write $(N,L,K,L^\prime,K^\prime)\Rightarrow \infty$. The macroscopic field, $\{S^{N,c}_v\}_{v\in V_N}$, is a centred Gaussian field with covariance matrix $\Sigma^c$, with entries given by
\begin{align}
	\Sigma^c_{u,v}\coloneqq \sigma^2(0) \mathbb{E}\left[\phi^{KL}_{w_i} \phi^{KL}_{w_j}\right],\quad \text{for } u\in B_{N/KL,i},\, v\in B_{N/KL,j},
\end{align}
where $\{\phi^{KL}_v\}_{v\in V_{KL}}$ is a DGFF on $V_{KL}$. It captures the macroscopic dependence.
The microscopic or ``bottom field``, $\{S^{N,b}_v\}_{v\in V_N}$, is a centred Gaussian field with covariance matrix $\Sigma^b$ defined entry-wise as
\begin{align}
	\Sigma^b_{u,v}\coloneqq 
	\begin{cases}
		\sigma^2(1) \mathbb{E}\left[\phi^{K^\prime L^\prime}_{u- v_{K^\prime L^\prime,i}} \phi^{K^\prime L^\prime}_{v-v_{K^\prime L^\prime,i}} \right],\quad &\text{if }u,v\in B_{K^\prime L^\prime,i}\\
		0, & \text{else},
	\end{cases}
\end{align}
where $\{\phi^{K^\prime L^\prime}_v \}_{v\in V_{K^\prime L^\prime}}$ is a DGFF on $V_{K^\prime L^\prime}$. It captures ``local'' correlations.
The third centred Gaussian field, $\{S^{N,m}_v \}_{v\in V_N}$, approximates the ``intermediate'' scales. It is a modified inhomogeneous branching random walk, defined pointwise as
\begin{align}
	S^{N,m}_v \coloneqq \sum_{j=k^\prime +l^\prime}^{n-l-k} \sum_{B \in \mathcal{B}_j(v_{K^\prime L^\prime,i^\prime})} 2^{-j}\sqrt{\log 2}b^N_{i,j,B} \int_{n-j-1}^{n-j}\sigma\left(\frac{s}{n}\right)\mathrm{d}s ,\quad \text{for } v\in B_{N/KL,i}\cap B_{K^\prime L^\prime, i^\prime},
\end{align}
with $\{b^N_{i,j,B}:  B\in \cup_{i^\prime} \mathcal{B}_j(v_{K^\prime L^\prime, i^\prime}), i=1,\dotsc, (KL)^2,\, j=1,\dotsc, (N/K^\prime L^\prime)^2,\ \}$ being a family of independent standard Gaussian random variables and where $\mathcal{B}_j(v_{K^\prime L^\prime, i^\prime})$ is the collection of boxes, $B\subset V_N$, of side length $2^j$ and lower left corner in $V_N$, that contain the element $v_{K^\prime L^\prime,i^\prime}$. In order to avoid boundary effects, we restrict our considerations onto a slightly smaller set, which is defined next. Consider the disjoint union of $N/L-$ and $L-$boxes, that is $\mathcal{B}_{N/L}=\{B_{N/L,i}:\, i=1\dotsc,L^2 \}$ and $\mathcal{B}_L=\{B_{L,i}:\, i=1,\dotsc,(N/L)^2\}$. Analogously, let $v_{N/L,i}$ and $v_{L,i}$ be the bottom left corners of the boxes $B_{N/L,i},\, B_{L,i}$ containing $v$. For a box $B$, let $B^\delta \subset B$ be the set $B^\delta=\{v\in B:\, \min_{z\in \partial B} \|v-z\|\geq \delta l_B \},$ where $l_B$ denotes the side length of the box $B$. Finally, let
\begin{align}
V_{N,\delta}^{*}\coloneqq\{\mathop{\cup}_{1\leq i \leq L^2}B^{\delta}_{N/L,i}\}\cap \{\mathop{\cup}_{1\leq i \leq (KL)^2} B^{\delta}_{N/KL,i} \}\cap \{\mathop{\cup}_{1\leq i \leq (N/L)^2}B^{\delta}_{L,i} \}\cap \{\mathop{\cup}_{1\leq i \leq (N/KL)^2}B^{\delta}_{KL,i}\}.
\end{align}
The next lemma ensures that the sum of the three fields, $\{S^{N,c}_v\}_{v\in V_N}, \{S^{N,m}_v\}_{v\in V_N}, \{S^{N,b}_v\}_{v\in V_N}$, approximates well the scale-inhomogeneous DGFF, $\{\psi^N_v\}_{v\in V_N}$.
\begin{lemma}\cite[Lemma~5.2, Lemma~5.3]{paper2}\label{lemma:4.1}
	There are non-negative uniformly bounded sequences of constants $a_{K^\prime L^\prime,\bar{v}}$ and a family of i.i.d. Gaussians $\{\Theta_j\}_{j=1,\dotsc, (N/K^\prime L^\prime)^2}$, such that, for $v\in B_{K^\prime L^\prime,j}$, $v\equiv \bar{v} \mod K^\prime L^\prime$, i.e. $\bar{v}=v-v_{K^\prime L^\prime,j}$, and when setting
	\begin{align}\label{equation:approximating_field}
	S^N_v \coloneqq S^{N,c}_v+ S^{N,m}_v+ S^{N,b}_v +a_{K^\prime, L^\prime, j} \Theta_j,
	\end{align}
	we have
	\begin{align}
	\limsup\limits_{(N,L,K,L^\prime,K^\prime)\Rightarrow \infty}\left| \mathrm{Var}\left(S^N_v \right) -\mathrm{Var}\left(\psi^N_v\right)-4\alpha \right|=0,
	\end{align}
	for some $\alpha>0$.
	Further, there exists a sequence $\{\epsilon_{N,KL,K^\prime L^\prime}^{'}\geq 0 \}$ with $\limsup\limits_{(N,L,K,L^\prime,K^\prime)\Rightarrow \infty} \epsilon_{N,KL,K^\prime L^\prime}^{'}=0$ and bounded constants $C_{\delta},C>0$, such that for all $u,v\in V_{N,\delta}^{*}:$
	\begin{enumerate}
		\item If $u,v\in B_{L^\prime,i},$ then $\left|\mathbbm{E}\left[\left(S^{N}_u-S^N_v\right)^2\right]- \mathbbm{E}\left[\left(\psi^N_u-\psi^N_v \right)^2\right]\right|\leq \epsilon_{N,KL,K^\prime L^\prime}^{'}.$
		\item If $u\in B_{N/L,i},$ $v\in B_{N/L,j}$ with $i\neq j$, then $\left|\mathbbm{E}\left[S^N_u S^N_v \right]-\mathbbm{E}\left[\psi^N_u\psi^N_v\right]\right|\leq \epsilon_{N,KL,K^\prime L^\prime}^{'}.$
		\item In all other cases, that is if $u,v\in B_{N/L,i}$ but $u\in B_{L^\prime,i^\prime}$ and $v\in B_{L^\prime,j^\prime}$ for some $i^\prime \neq j^\prime$, it holds that $\left|\mathbbm{E}\left[S^N_uS^N_v \right]-\mathbbm{E}\left[\psi^N_u\psi^N_v \right]\right|\leq C_{\delta}+ 40\alpha$.
	\end{enumerate}
\end{lemma}
The field, $\{S^N_v\}_{v\in V_N}$, defined in \eqref{equation:approximating_field} is the approximating $3-$field we work with.
	\subsection{Reduction to approximating field}
	In the following, we generalize the approximation results from \cite{paper2} to the case of countably many local maxima. We show that the local maxima of $\{\psi^N_v\}_{v\in V_N}$ are well approximated by those of $\{S^N_v\}_{v\in V_N}$. As we need to compare probability measures on $\mathbb{R}^p$, we use the L\'{e}vy-Prokhorov metric $d(\cdot,\cdot)$, to measure distances between probability measures on $\mathbb{R}^p$. For two probability measures, $\mu$ and $ \nu$, it is given by
\begin{align}
	d(\mu,\nu)= \inf \{\delta>0:\,\mu(B)\leq \nu(B^\delta) +\delta \text{ for all open sets } B \},
\end{align}
where $B^\delta=\{ y\in \mathbb{R}^p: |x-y|<\delta, \text{ for some } x\in B\}$. Further, let
\begin{align}
	\tilde{d}(\mu, \nu)=\inf  \{\delta>0:\, \mu((x_1,\infty),\dotsc,(x_p,\infty))\leq \nu((x_1-\delta,\infty),\dotsc, (x_p-\delta,\infty)) +\delta,\forall (x_1,\dotsc, x_p)\in \mathbb{R}^p \},
\end{align}
which is a measure for stochastic domination.
In particular, if $\tilde{d}(\mu,\nu)=0$, then $\nu$ stochastically dominates $\mu$. Note, unlike $d(\cdot,\cdot)$, $\tilde{d}(\cdot,\cdot)$ is not symmetric. Abusing notation, we write for random vector $X,Y$ with laws $\mu_X, \mu_Y$, $d(X,Y)$ instead of $d(\mu_X,\mu_Y)$ and likewise for $\tilde{d}$. Fix $r\in \mathbb{N}$ and let $\mathcal{B}_r$ of $V_{\lfloor N/r\rfloor r}$ into sub-boxes of side length $r$. Let $\mathcal{B}= \cup_{r\in \mathbb{N},r\leq N}\mathcal{B}_r$ and $\{g_b\}_{B\in \mathcal{B}}$ be a collection of i.i.d. standard Gaussian random variables. For $v\in V_N$, denote by $B_r(v)\in \mathcal{B}_r$ the box containing $v$. For $r_1,r_2\in \mathbb{N}$, $r_1,r_2\leq N$, $A\subset [0,1]^2$, $s_1.s_2\in \mathbb{R}_{+}$, we write
\begin{align}
	\bar{\psi}^{*}_{N,A}\coloneqq \max\limits_{v\in V_N: v/N\in A}\psi^N_v+ s_1 g_{B_{v,r_1}}+s_2 g_{B_{v,N/r_2}},
\end{align}
and for a general field $\{g^N_v\}_{v\in V_N}$,
\begin{align}
	g^{*}_{N,A}\coloneqq \max_{v\in V_N:v/N\in A} g^N_v.
\end{align}
Fix $p\in \mathbb{N}$ and disjoint, open, non-empty, simply connected sets $A_1,\dotsc, A_p\subset[0,1]^2$.
\begin{lemma}\label{lemma:approx_limiting_laws}
	For $s=(s_1,s_2)\in \mathbb{R}_{+}^2$, it holds
	\begin{align}
		\limsup\limits_{r_1,r_2\rightarrow \infty}\limsup\limits_{N\rightarrow \infty} d((\psi^{*}_{N,A_i}-m_N)_{1\leq i \leq p},(\bar{\psi}^{*}_{N,A_i}-m_N -\|s\|^2_2)_{1\leq i \leq p})=0.
	\end{align}
\end{lemma}
For the proof of \autoref{lemma:approx_limiting_laws} we need some additional estimates.
\begin{lemma}\label{lemma:generalized_perturbation_of_max_lemma}
	Let $\{\bar{\psi}^N_v\}_{v\in V_N}$ be a centred Gaussian field and $c>0$ a constant, such that, for any $v,w\in V_N$, $\left| \mathbb{E}\left[\bar{\psi}^N_v \bar{\psi}^N_w\right]-\mathbb{E}\left[\psi^N_v \psi^N_w\right]\right| \leq c.$ Moreover, let $A\subset [0,1]^2$ be an open, non-empty subset and $\{g^N_v \}_{v\in V_N}$ be a collection of independent random variables, such that
	\begin{align}\label{equation:assumption_gaussian_tail_rv}
		\mathbb{P}\left(g^N_v \geq 1+y\right)\leq e^{-y^2}\quad \text{for } v\in V_N.
	\end{align}
	Then, there is a constant $C=C(\alpha)>0$ such that, for any $\epsilon>0$, $N\in \mathbb{N}$ and $x\geq -\epsilon^{1/2}$,
	\begin{align}
		\mathbb{P}\left(\max_{v\in V_N: v/N \in  A} (\bar{\psi}^N_v +\epsilon g^N_v)\geq m_N +x \right)\leq \mathbb{P}\left( \max_{v\in V_N:v/N\in A} \bar{\psi}^N_v \geq m_N+x-\sqrt{\epsilon}\right)(Ce^{-C^{-1}\epsilon^{-1}}).
	\end{align}
	\begin{proof}
		Set $\Gamma_y \coloneqq \{ v\in V_N:v/N\in A, y/2 \leq \epsilon g^N_v \leq y\}.$ Then,
		\begin{align}\label{equation:3.14}
			\mathbb{P}\left(\max_{v\in V_N:v/N \in A}(\bar{\psi}^N_v +\epsilon g^N_v)\geq m_N +x \right)&\leq \mathbb{P}\left(\max_{v\in V_N:v/N\in A}\bar{\psi}^N_v \geq m_N+x- \sqrt{\epsilon}\right)\nonumber\\&\quad
			+ \sum_{i=0}^{\infty}\mathbb{E}\left[ \mathbb{P}\left(\max_{v\in \Gamma_{2^i \sqrt{\epsilon}}}\bar{\psi}^N_v \geq m_N +x -2^i \sqrt{\epsilon} | \Gamma_{2^i \sqrt{\epsilon}}\right) \right].
		\end{align}
		By \cite[Proposition~5.1]{paper2}, the second term on the right hand side in \eqref{equation:3.14} is bounded from above by
		\begin{align}\label{equation:sum_prob_1}
			\sum_{i=0}^{\infty} \mathbb{E}\left[\mathbb{P}\left(\max_{v\in V_N:v/N\in A}\bar{\psi}^N_v \geq m_N +x -2^i \sqrt{\epsilon} | \Gamma_{2^i \sqrt{\epsilon}} \right) \right]\leq \tilde{c}e^{-cx} \sum_{i=0}^{\infty} \mathbb{E}\left[| \Gamma_{2^i \sqrt{\epsilon}}|/ | \{v\in V_N:v/N\in A\}| \right] e^{c 2^i \sqrt{\epsilon}},
		\end{align} 
		where $\tilde{c}>0$ is a finite constant. By assumption \eqref{equation:assumption_gaussian_tail_rv}, one has
		\begin{align}\label{equation:3.16}
			\mathbb{E}\left[| \Gamma_{2^i\sqrt{\epsilon}}|/ |\{v\in V_N : v/N \in A\}| \right]\leq e^{-4^i (C\epsilon)^{-1}}.
		\end{align}
		Thus, using  \eqref{equation:3.16}, \eqref{equation:sum_prob_1} is bounded from above by
		\begin{align}
			\tilde{c}e^{-cx} e^{-(C\epsilon)^{-1}}.
		\end{align}
		This concludes the proof of \autoref{lemma:generalized_perturbation_of_max_lemma}.
	\end{proof}
\end{lemma}
\begin{lemma}\label{lemma:stochastic_domination_cov_comp}
	Let $\{\bar{\psi}^N_v\}_{v\in V_N}$ be a centred Gaussian field satisfying
	\begin{align}\label{equation:3.19}
		| \mathrm{Var}\, \psi^N_v - \mathrm{Var}\, \tilde{\psi}^N_v | \leq \epsilon.
	\end{align}
	Further, fix some $p\in \mathbb{N}$, and disjoint open, non-empty sets $A_1,\dotsc, A_p\subset [0,1]^2$.
If
\begin{align}\label{equation:3.20}
	\mathbb{E}\left[ \tilde{\psi}^N_v \tilde{\psi}^N_w\right]\leq \mathbb{E}\left[\psi^N_v \psi^N_w \right] + \epsilon,
\end{align}
then
\begin{align}\label{equation:4.20}
	\limsup\limits_{N\rightarrow \infty} \tilde{d}\left((\psi^{*}_{N, A_1}-m_N,\dotsc,\psi^{*}_{N,A_p}-m_N), (\tilde{\psi}^{*}_{N,A_1}-m_N,\dotsc, \tilde{\psi}^{*}_{N,A_p}-m_N) \right)\leq l(\epsilon),
\end{align}
and else if,
\begin{align}
\mathbb{E}\left[ \tilde{\psi}^N_v \tilde{\psi}^N_w\right]+\epsilon\geq \mathbb{E}\left[\psi^N_v \psi^N_w \right],
\end{align}
 then
 \begin{align}\label{equation:4.22}
 \limsup\limits_{N\rightarrow \infty} \tilde{d}\left((\tilde{\psi}^{*}_{N, A_1}-m_N,\dotsc,\tilde{\psi}^{*}_{N,A_p}-m_N), (\psi^{*}_{N,A_1}-m_N,\dotsc, \psi^{*}_{N,A_p}-m_N) \right)\leq l(\epsilon),
 \end{align}
where $l(\epsilon)\rightarrow 0$ as $\epsilon \rightarrow 0$.
\begin{proof}
	Let $\{\psi^N_v \}_{v\in V_N}, \, \{\tilde{\psi}^N_v \}_{v\in V_N}$ satisfy relations \eqref{equation:3.19} and \eqref{equation:3.20}. Let $\Phi, \{\Phi^N_v\}_{v\in V_N}$ two independent standard Gaussian random variables, and $\epsilon^{*}(\epsilon)>0$. For $v\in V_N$, set
	\begin{align}
		\psi^{N,lw,\epsilon^{*}}_v&=\left(1-\frac{\epsilon^{*}}{\log N}\right)\psi^N_v + \epsilon^{N,\prime} \Phi,\\
		\tilde{\psi}^{N,up,\epsilon^{*}}_v&=\left(1-\frac{\epsilon^{*}}{\log N}\right)\tilde{\psi}^N_v + \epsilon^{N,\prime \prime}_v \Phi^N_v,
	\end{align}
	where we can choose, as in the proof of \cite[Lemma~5.6]{paper2}, $\epsilon^{*}$, $\epsilon^{N,\prime}_v=\epsilon^{N,\prime}_v(\epsilon, \epsilon^{*})$ and $\epsilon^{N,\prime \prime}_v=\epsilon^{N,\prime \prime}_v(\epsilon, \epsilon^{*})$ all non-negative and tending to $0$ as $\epsilon\rightarrow 0$, such that
	\begin{align}
		\mathrm{Var}\left[\psi^{N,lw,\epsilon^{*}}_v\right]=\mathrm{Var}\left[\tilde{\psi}^{N,up,\epsilon^{*}}_v\right]=\mathrm{Var}\left[\psi^N_v\right] +\epsilon, \quad\forall v\in V_N
	\end{align}
	and
	\begin{align}
		\mathbb{E}\left[\psi^{N,lw,\epsilon^{*}}_v \psi^{N,lw,\epsilon^{*}}_w\right]\geq \mathbb{E}\left[\tilde{\psi}^{N,up,\epsilon^{*}}_v \tilde{\psi}^{N,up,\epsilon^{*}}_w\right], \quad \forall v,w\in V_N.
	\end{align}
	An application of Slepian's lemma for vectors (\autoref{theorem:vector_Slepian}), gives
	\begin{align}\label{equation:6.25}
		\tilde{d}\left( (\psi^{*}_{N,lw,\epsilon^{*},A_1}-m_N,\dotsc, \psi^{*}_{N,lw,\epsilon^{*},A_p}-m_N),(\tilde{\psi}^{*}_{N,up,\epsilon^{*},A_1}-m_N,\dotsc, \tilde{\psi}^{*}_{N,up,\epsilon^{*},A_1}-m_N) \right)=0.
	\end{align}
	By \autoref{lemma:generalized_perturbation_of_max_lemma}, we obtain, for $x_1,\dotsc, x_p \in \mathbb{R}$,
	\begin{align}\label{equation:6.26}
		\mathbb{P}\left(\tilde{\psi}^{*}_{N,up,\epsilon^{*},A_i}-m_N\geq x_i, \, 1\leq i \leq p \right)&\leq \mathbb{P}\left(\psi^{*}_{N,A_i}-m_N\geq x_i-\sqrt{\max_{w\in V_N}\epsilon^{N,\prime \prime}_w},\, 1\leq i \leq p \right)\nonumber\\ &\qquad\times Ce^{-(C \max_{w\in V_N}\epsilon^{N,\prime \prime}_w )^{-1}}.
	\end{align}
	Since $\lim\limits_{\epsilon \rightarrow 0}\max_{w\in V_N}\epsilon^{N,\prime \prime}_w=0$ this implies \eqref{equation:4.20}. \eqref{equation:4.22} can be proved the same way by switching the roles of $\{\psi^N_v\}_{v\in V_N}$ and $\{\tilde{\psi}^N_v \}_{v\in V_N}$. We omit further details.
\end{proof}
\end{lemma}
\begin{prop}\label{proposition:3.5}
	Let $\tilde{\sigma}\in (0,\infty)^2$, $r=(r_1,r_2)\in (0,\infty)^2$, and $\{\psi^{N,r, \tilde{\sigma}}_v:\, v\in V_N \}$ as well as $\{\psi^{N,\tilde{\sigma},*}_v: \, v\in V_N  \}$ be two Gaussian fields given by
	\begin{align}
		\psi^{N,r,\tilde{\sigma}}_v=\psi^N_v + \tilde{\sigma}_1 g_{B_{v,r_1}}+ \tilde{\sigma}_2 g_{B_{v,N/r_2}}, \quad \text{for } v\in V_n
	\end{align}
	and
	\begin{align}\label{equation:3.31}
		\psi^{N,\tilde{\sigma},*}_v=\psi^N_v +\sqrt{\frac{\|\tilde{\sigma}\|^2_2}{\log(N)}}\tilde{\psi}^N_v,\quad \text{for } v\in V_N
	\end{align}
	where $\{\psi^N_v\}_{v\in V_N},\,\{\tilde{\psi}^N_v\}_{v\in V_N}$ are two independent scale-inhomogeneous DGFFs, satisfying \autoref{assumption:1}, and where $\{g_B\}_{B\in \mathcal{B}}$ is a collection of independent standard Gaussians. For a set $A\subset [0,1]^2$, we write $M_{N,A,r_1,r_2,\tilde{\sigma}}=\underset{v \in V_N:v/N\in A}{\max} \psi^{N,r,\tilde{\sigma}}_v$ and likewise, $M_{N,A,\tilde{\sigma},*}=\underset{v \in V_N:v/N\in A}{\max}\psi^{N,\tilde{\sigma},*}_v$.
	Then, for any $p\in \mathbb{N}$, and any collection of disjoint, open and non-empty $A_1,\dotsc, A_p \subset[0,1]^2$,
	\begin{align}
		\limsup\limits_{N\rightarrow \infty} d\left((M_{N,A_1,r_1,r_2,\tilde{\sigma}}-m_N,\dotsc, M_{N,A_p,r_1,r_2,\tilde{\sigma}}-m_N) ,(M_{N,A_1,\tilde{\sigma},*}-m_N,\dotsc, M_{N,A_p,\tilde{\sigma}*}-m_N)\right) =0,
	\end{align}
	as $r_1,r_2\rightarrow \infty$.
	\begin{proof}
		The proof is a straightforward adaptation of the proof of \cite[Proposition~B.2]{paper2}.
		Decompose $V_N$ into boxes $B$ of side length $N/r_2$ and call their collection $\mathcal{B}$. Further, for $\delta \in (0,1)$ and $B\in \mathcal{B}$, let $B_{\delta}$ be the box with the identical centre as $B$, and reduced side length $(1-\delta)N/r_2$. Then, we set $V_{N,\delta}=\cup_{B\in \mathcal{B}} B_{\delta}$. The corresponding maxima over are called $M_{N,A,r_1,r_2,\tilde{\sigma},\delta}= \underset{v\in V_{N,\delta}: v/N\in A}{\max}\psi^{N,r,\tilde{\sigma}}_v$ and $M_{N,A,\tilde{\sigma},*}=\underset{v\in V_{N,\delta}: v/N\in A}{\max}\psi^{N,\tilde{\sigma},*}_v$.
		\cite[Proposition~5.1]{paper2} shows that it suffices to consider the maxima on the slightly smaller sets, i.e. one has
		\begin{align}
			\lim\limits_{\delta\rightarrow 0} \lim\limits_{N\rightarrow \infty}\mathbb{P}\left(M_{N,A_1,r_1,r_2,\tilde{\sigma},\delta}\neq M_{N,A_1,r_1,r_2,\tilde{\sigma}},\dotsc,M_{N,A_p,r_1,r_2,\tilde{\sigma},\delta}\neq M_{N,A_p,r_1,r_2,\tilde{\sigma}} \right)\nonumber\\
			= \lim\limits_{\delta \rightarrow 0}\lim\limits_{N\rightarrow \infty}\mathbb{P}\left(M_{N,A_1,\tilde{\sigma},*,\delta}\neq M_{N,A_1,\tilde{\sigma},*},\dotsc, M_{N,A_p,\tilde{\sigma},*,\delta}\neq M_{N,A_p,\tilde{\sigma},*} \right)=0.
		\end{align}
		Next, we claim that the maximum is essentially determined by the maximum of the unperturbed scale-inhomogeneous DGFF, $\{\psi^N_v\}_{v\in V_N}$. For $B\in \mathcal{B}$, let $z_B$ be the unique element, such that
		\begin{align}
			\psi^N_{z_B}=\max_{v \in B_{\delta}} \psi^N_v.
		\end{align}
		The claim is that
		\begin{align}\label{equation:6.21}
			\lim\limits_{r_1,r_2\rightarrow \infty}\lim\limits_{N\rightarrow \infty}&\mathbb{P}\left(|M_{N,A_i,r_1,r_2,\tilde{\sigma},\delta}-\max_{B\in \mathcal{B},B\subset N A_i}\psi^{N,r,\tilde{\sigma}}_{z_B}|\geq \frac{1}{\log n}:\,1\leq i\leq p \right)\nonumber\\
			&=\limsup\limits_{N\rightarrow \infty} \mathbb{P}\left(|M_{N,A_i,\tilde{\sigma},*,\delta}-\max_{B\in \mathcal{B},B\subset N A_i}\psi^{N,\tilde{\sigma},*}_{z_B}|\geq \frac{1}{\log n}:\, 1\leq i \leq p\right)=0.
		\end{align}
		In the following, we show that none of the events in the probabilities in \eqref{equation:6.21} can occur. It suffices to show that none of the following events can happen. For $i\in \{1,\dotsc,p\}$, let
		\begin{align}
			E_1^{(i)}=&\{M_{N,A_i,r_1,r_2,\tilde{\sigma},\delta } \notin (m_N-C,m_N+C) \} \cup \{M_{N,A_i,\tilde{\sigma},*,\delta}\notin (m_N-C,m_N+C) \} \\
			E_2^{(i)}=&\{\exists u,v\in V_N: u,v/N\in A_i, \|u-v\| \in (r,N/r)\text{ and } \min(\psi^N_u,\psi^N_v)>m_N -c\log n \} \\
			E_3^{(i)}=&\tilde{E}_3^{(i)}\cup \bar{E}_3^{(i)},\text{ where } \tilde{E}^{(i)}_3=\{\omega:\exists v\in V_N,\, v/N\in A_i: \psi^{N,r,\tilde{\sigma}}_v=M_{N,A_i,r_1,r_2,\tilde{\sigma},\delta },\, \psi^N_v\leq m_N -c\log n \},\nonumber\\ & \bar{E}^{(i)}_3=\{\omega:\exists v\in V_N,\, v/N\in A_i: \psi^{N,\tilde{\sigma},*}_v=M_{N,A_i,\tilde{\sigma},*,\delta},\, \psi^N_v\leq m_N -c\log n \} \\
			E_4^{(i)}=&\left\{\exists v\in B\in \mathcal{B}\subset N A_i: \psi^N_v\geq m_N -c\log n \text{ and }
			\sqrt{\frac{\|\tilde{\sigma}\|^2_2}{\log N}} \tilde{\psi}^N_v- \sqrt{\frac{\| \tilde{\sigma}\|^2_2}{\log N}} \tilde{\psi}^N_{z_B}\geq 1/\log n \right\} .
		\end{align}
		The events $E_2,E_3$ and $E_4$ in the proof of \cite[Proposition~B.2]{paper2} include the corresponding events, $E^{(i)}_2, E^{(i)}_3, E^{(i)}_4$, we are considering here, and so we know that the probability of their occurrence tends to zero. So, we are left with bounding the events $E^{(i)}_1$. First note that it suffices to consider the scale-inhomogeneous DGFF, as the other terms are centred Gaussians with uniformly bounded variance. Since maximizing over a subset, we have, for any $i\in \{1,\dotsc, p\}$,
		\begin{align}\label{equation:3.40}
			\mathbb{P}\left(\max_{v\in V_N:\, v/N\in A_i} \psi^N_v > m_N +C \right)\leq \mathbb{P}\left(\max_{v\in V_N}\psi^N_v >m_N+C \right).
		\end{align}
		By tightness of the centred maximum \cite[(2.2)]{paper2}, \eqref{equation:3.40} tends to $0$ as $C\rightarrow \infty$, uniformly in $N$. Hence to show \eqref{equation:6.21}, it suffices to prove, for any $i\in \{1,\dotsc, p\}$,
		\begin{align}\label{equation:6.27}
			\lim\limits_{C\rightarrow \infty}\lim\limits_{N\rightarrow \infty} \mathbb{P}\left( \max_{v\in V_N: v/N \in A_i} \psi^N_v \leq m_N -C\right)=0.
		\end{align}
		Assume otherwise, then there is a subsequence $\{N_k\}_{k\in \mathbb{N}}$, a sequence $C_N\rightarrow \infty$ as $N\rightarrow \infty$ and a constant $\epsilon>0$, such that, for any $k\in \mathbb{N}$,
		\begin{align}\label{equation:abc}
		\mathbb{P}\left(\max_{v\in V_{N_k}:\, v/N_k \in A_i} \psi^{N_k}_v \leq m_{N_k}- C_{N_k} \right)\geq \epsilon.
		\end{align} 
		We can further assume that $A_i \subset [0,1]^2$ is a box, otherwise pick the largest box that fits into $A_i$.
		We can decompose $[0,1]^2$ into disjoint translations of $A^{(j)}_i$, that we possible need to cut with $[0,1]^2$. For each $A^{(j)}_i N$ we consider an independent copy of $\{\psi^N_v\}_{v\in V_N}$, called $\{\psi^{N,j}_v\}_{v\in V_N}$. By translation invariance, for each of these \eqref{equation:abc} holds.
		By Gaussian comparison, independence and \eqref{equation:abc}, we have
		\begin{align}\label{equation:ab}
			\mathbb{P}\left(\max_{v\in V_{N_k}}\psi^{N_k}_v\leq m_{N_k}-C_{N_k}\right)\geq \mathbb{P}\left(\max_j \max_{v\in A^{(j)}_i N_k}\psi^{N_k,j}_v\leq m_{N_k}-C_{N_k}\right)>0.
		\end{align}
		By tightness of $\{\max_{v\in V_N}\psi^N_v-m_N \}_{N\in \mathbb{N}}$, the left-hand side of \eqref{equation:ab} tends to zero, which is a contradiction. Thus, this yields \eqref{equation:6.27}, which concludes the proof of \autoref{proposition:3.5}.
	\end{proof}
\end{prop}
\autoref{lemma:stochastic_domination_cov_comp} and \autoref{proposition:3.5} allow us to prove \autoref{lemma:approx_limiting_laws}.
\begin{proof}[Proof of \autoref{lemma:approx_limiting_laws}:]
	Define for $v\in V_N$, $\bar{\psi}^{N,\tilde{\sigma}}_v=\left(1+ \frac{\|\tilde{\sigma}\|^2}{\log(N)}\right)\psi^N_v$, and for $A\subset [0,1]^2$ open and non-empty,  $\bar{M}_{N,A,\tilde{\sigma}}= \max_{v\in V_N:v/N \in A}\bar{\psi}^{N,\tilde{\sigma}}_v$ and set $M_{N,A}=\max_{v \in V_N:v/N\in A}\psi^N_v$. \eqref{equation:3.40} together with tightness of the centred maximum \cite[(2.2)]{paper2} and \eqref{equation:6.27} implies,
	\begin{align}
		\mathbb{E}\left[\bar{M}_{N,A_i,\tilde{\sigma}} \right]= \mathbb{E}\left[M_{N,A_i}\right]+2 \|\tilde{\sigma}\|^2_2 +o(1),
	\end{align}
	and
	\begin{align}
		\lim\limits_{N\rightarrow \infty} d(M_{N,A_i}- \mathbb{E}\left[M_{N,A_i}\right], \bar{M}_{N,A_i,\tilde{\sigma}}-\mathbb{E}\left[\bar{M}_{N,A_i,\tilde{\sigma}}\right])=0.
	\end{align}
	Next, we consider the field, $\{\psi^{N,\tilde{\sigma},*}_v\}_{v\in V_N}$, defined in \eqref{equation:3.31}. For $i\in \{1,\dotsc,p\}$, set $M_{N,A_i,\tilde{\sigma},*}=\max_{v \in V_N:\, v/N \in A_i}\psi^{N,\tilde{\sigma},*}_v$. In distribution, $\{\psi^{N,\tilde{\sigma},*}_{v}\}_{v\in V_N}$ can be written as a sum of $\{\bar{\psi}^{N,\tilde{\sigma}}_v\}_{v\in V_N}$ and an independent centred Gaussian field with variances of order $O((1/\log N)^3)$. Thus, by Gaussian comparison,
	\begin{align}
		\mathbb{E}\left[\bar{M}_{N,A_i,\tilde{\sigma}}\right]=\mathbb{E}\left[M_{N,A_i,\tilde{\sigma},*}\right]+o(1)
	\end{align}
	and
	\begin{align}\label{equation:3.46}
		\lim\limits_{N\rightarrow \infty}d\left(\left(\bar{M}_{N,A_i,\tilde{\sigma}}-\mathbb{E}\left[\bar{M}_{N,A_i,\tilde{\sigma}}\right]\right)_{1\leq i\leq p},\left(\bar{M}_{N,A_i,\tilde{\sigma},*}-\mathbb{E}\left[\bar{M}_{N,A_i,\tilde{\sigma},*}\right]\right)_{1\leq i\leq p}\right)=0.
	\end{align}
	Combining \eqref{equation:3.46} with \autoref{proposition:3.5} and applying the triangle inequality, one concludes the proof of \autoref{lemma:approx_limiting_laws}.
\end{proof}
Finally, we are able to deduce the key result in this subsection.
\begin{lemma}\label{lemma:reduction_multiple_to_3field}
	Let $p\in \mathbb{N}$, and $A_1,\dotsc,A_p\subset [0,1]^2$ be  disjoint, open and non-empty. Then,
	\begin{align}
	\limsup\limits_{(N,L,K,L^\prime,K^\prime)\Rightarrow \infty} d\left((\psi^{*}_{N,A_i}-m_N)_{ 1\leq i \leq p},(S^{*}_{N,A_i}-m_N -4\alpha)_{1\leq i \leq p} \right)=0.
	\end{align}
	\begin{proof}
		We refrain from giving the proof, as it follows in complete analogy to \cite[Lemma~5.4]{paper2}. Instead of using \cite[Lemma~5.6]{paper2} in the proof, one replaces it by its multi-dimensional analogue, \autoref{lemma:stochastic_domination_cov_comp}.
	\end{proof}
\end{lemma}
This reduces the proof of convergence in law of multiple local maxima of the scale-inhomogeneous DGFF to the structurally simpler field, $\{S^N_v\}_{v\in V_N}$, as it decouples microscopic and macroscopic dependence.
	\subsection{Coupling to independent random variables}
	Recall $\underline{A}=(A_1,\dotsc,A_p)$ is a collection of disjoint open, non-empty, simply-connected subsets of $[0,1]^2$, for some fixed $p\in \mathbb{N}$.
Further, we have tiled $V_N$ with boxes $B_{N/KL,i}$ of side length $N/KL$. Instead of considering the maximum over the sets $\{v\in V_N: v/N\in A_i \}$, we want to work with the $B_{N/KL}$-boxes. Thus, for any $i\in \{1,\dotsc, p\}$, let $T^{(KL)}_i \subset \{1,\dotsc, (KL)^2\}$ denote the maximal index set, such that $j\in T^{(KL)}_i$ implies $B_{N/KL,j}/N \subset A_i$, i.e.
\begin{align}
	\cup_{j\in T^{(KL)}_i} B_{N/KL,j}/N \subset A_i.
\end{align}
Further, it is immediate to see that for all $1\leq i \leq p$
\begin{align}\label{equation:4.48}
	\frac{|N A_i \setminus \cup_{j\in T^{(KL)}_i} B_{N/KL,j} |}{|N A_i |}\rightarrow 0,
\end{align}
as we let $N, K,L$ tend to infinity in this order. In particular,
\begin{align}\label{equation:4.49}
	\mathbb{P}\left(\max_{v\in \cup_i{i=1}^p \left( A_i\setminus \cup_{j\in T^{(KL)}_i} B_{N/KL,j}\right) }\psi^N_v\geq m_N+z\right)&\leq \sum_{i=1}^{p} |N A_i \setminus \cup_{j\in T^{(KL)}_i} B_{N/KL,j} | \sup_{v\in V_N} \mathbb{P}\left(\psi^N_v\geq m_N+z\right)\nonumber\\
	&\leq C \sum_{i=1}^p \frac{|N A_i \setminus \cup_{j\in T^{(KL)}_i} B_{N/KL,j} |}{N^2}e^{-2z},
\end{align}
which, by \eqref{equation:4.48}, converges to zero as $N\rightarrow \infty$.
Next, we construct random variables that do not depend on $N$ and that we couple to the local maxima of $\{S^N_v\}_{v\in V_N}$ on $\cup_{j\in T^{(KL)}_1}B_{N/KL,j},\dotsc, \cup_{j\in T^{(KL)}_p}B_{N/KL,j}$.
We set $A_i^\prime\coloneqq\cup_{j\in T^{(KL)}_i} B_{N/KL,j}$, and $S^{N,f}_v\coloneqq S^N_v-S^{N,c}v,$ for $v\in V_N$.
Let $\{\varrho_{R,i}:\, 1\leq i\leq R \}$ be a collection of independent Bernoulli random variables with
\begin{align}\label{equation:4.50}
	\mathbb{P}\left(\varrho_{R,i}=1\right)= \beta^{*}_{K^\prime, L^\prime} e^{2 \bar{k}^{\gamma}} e^{2\bar{k}(\sigma^2(0)-1)},
\end{align}
where, by using \cite[Proposition~5.8]{paper2}, the constants $\beta^{*}_{K^\prime, L^\prime}$ are such that they satisfy,
\begin{align}
	\lim\limits_{z\rightarrow \infty}\limsup\limits_{(L^\prime,K^\prime,N)\Rightarrow \infty}\left|e^{2\log(2)\bar{k}(1-\sigma^2(0))}e^{-2\bar{k}^\gamma}e^{2z}\mathbb{P}\left(\max_{v\in B_{N/KL,i}}S^{N,f}_v\geq m_N(\bar{k},n)-\bar{k}^\gamma+z\right)-\beta^{*}_{K^\prime, L^\prime}\right|=0.
\end{align}
Moreover, there are constans $c_\alpha,C_\alpha>0$ such that $c_\alpha \leq \beta^{*}_{K^\prime, L^\prime}\leq C_\alpha$, where $\alpha$ is as in \autoref{lemma:4.1}, and the collection $\{\beta^{*}_{K^\prime, L^\prime}\}_{K^\prime, L^\prime \geq 0}$ depends on the variance only through $\sigma(1)$.
In addition, we specify an independent family of exponential random variables, $\{Y_{R,i}:\, 1\leq i\leq R\}$,
\begin{align}\label{equation:4.51}
	\mathbb{P}\left(Y_{R,i}\geq x \right)=e^{-2x}e^{2\bar{k}^{\gamma}},\quad \text{for } x\geq -\bar{k}^{\gamma}.
\end{align}
Also, let $\{Z_{R,i}\}_{1\leq i\leq R\ }$ be a centered Gaussian field with correlation kernel $\Sigma^c$.
For each $i\in \{1,\dotsc,p \}$, set
\begin{align}\label{equation:4.52}
	G^{(i)}_{L,K,L^\prime, K^\prime}\coloneqq \max_{\substack{j \in T^{(KL)}_i \\ \varrho_{R,j}=1}} (Y_{R,j}+ 2\log(KL)(1-\sigma^2(0))) + (Z_{R,j}-2\log(KL)).
\end{align}
We collect these in the vector
\begin{align}\label{equation:3.53}
	G^{*}_{\underline{A},L,K,L^\prime, K^\prime}\coloneqq\left(G^{(1)}_{L,K,L^\prime, K^\prime},\dotsc, G^{(p)}_{L,K,L^\prime, K^\prime} \right).
\end{align}
We denote the law of the random vector defined in \eqref{equation:3.53} by $\bar{\mu}_{L,K,L^\prime,K^\prime,\underline{A}}$, which does not depend on $N$.
Next, we show that $\bar{\mu}_{L,K,L^\prime,K^\prime,\underline{A}}$ converges to the same limit as $\mu_{N,\underline{A}}$, the law of
\begin{align}
	\left(\max_{v\in A_1^\prime}S^N_v-m_N, \dotsc, \max_{v\in A_p^\prime}S^N_v-m_N  \right).
\end{align}
Set $m_N(k,t)\coloneqq 2\log N \mathcal{I}_{\sigma^2}\left(\frac{k}{n},\frac{t}{n}\right)-\frac{(t\wedge(n-\bar{l}))\log n}{4(n-\bar{l})}$, for $k\leq n$ and $t\in [k,n]$.
\begin{thm}\label{theorem:coupling_proof}
	It holds that
	\begin{align}
		\limsup\limits_{(N,L,K,L^\prime, K^\prime)\Rightarrow \infty} d(\mu_{N,\underline{A}}, \bar{\mu}_{L,K,L^\prime,K^\prime,\underline{A}} )=0.
	\end{align}
	In particular, there exists $\mu_{\infty, \underline{A}}$ such that $\lim\limits_{N\rightarrow \infty} d(\mu_{N,\underline{A}}, \mu_{\infty, \underline{A}})=0.$
	\begin{proof}
		We follow the proof of \cite[Theorem~5.9]{paper2} that deals with the global maximum. Denote by $\tau^\prime_i=\arg \max_{v\in B_{N/KL,i}}S^N_v$, the a.s. unique point where the local maximum is achieved.
		By \cite[(5.50)]{paper2}, we have, for $1\leq i \leq p$,
		\begin{align}
			\limsup\limits_{(N,L,K,L^\prime, K^\prime)\Rightarrow \infty} \mathbb{P}\left(S^{N,f}_{\tau^\prime_i}\geq m_N(\bar{k},n)- \bar{k}^\gamma \right)=1.
		\end{align}
		Moreover, we know that the fine field values cannot be too large, i.e. let 
		\begin{align}
			\mathcal{E}=\cup_{1\leq i \leq R}\{\max_{v \in B_{N/KL,i}}S^{N,f}_v\geq m_N(\bar{k},n)+KL +\bar{k}^\gamma \}, \text{ and } \mathcal{E}^\prime=\cup_{1\leq i \leq R}\{Y_{R,i}\geq KL +\bar{k}^\gamma \}.
		\end{align}
		By \cite[(5.51)]{paper2} respectively \cite[(5.53)]{paper2}, we deduce
		\begin{align}
			\limsup\limits_{(N,L,K,L^\prime, K^\prime)\Rightarrow \infty} \mathbb{P}\left(\mathcal{E}\right)=0\text{ and } \limsup\limits_{(N,L,K,L^\prime, K^\prime)\Rightarrow \infty}\mathbb{P}\left(\mathcal{E}^\prime\right)=0.
		\end{align}
		This allows to couple the centred fine field, $\tilde{M}^f_{N,i}=\max_{v \in B_{N/KL,i}}S^{N,f}_i-m_N(\bar{k},n)$, to the approximating process $G^{(i)}_{L,K,L^\prime, K^\prime}$, defined in \eqref{equation:4.52}.
	 	By \cite[Proposition~5.8]{paper2}, there are $\epsilon^{*}_{N,KL,K^\prime L^\prime}>0$ with
		\begin{align}
			\limsup\limits_{(N,L,K,L^\prime, K^\prime)\Rightarrow \infty}  \epsilon^{*}_{N,KL,K^\prime L^\prime}=0,
		\end{align}
		such that, for some $|\overset{\diamond}{\epsilon}|\leq \epsilon^{*}_{N,KL,K^\prime L^\prime}/4$,
		\begin{align}
			\mathbb{P}\left(-\bar{k}^\gamma +\overset{\diamond}{\epsilon}\leq \tilde{M}^f_{N,i}\leq KL + \bar{k}^\gamma \right)= \mathbb{P}\left(\varrho_{R,i}=1,\, Y_{R,i}\leq KL +\bar{k}^\gamma\right),
		\end{align}
		and such that for all $t$ with $-\bar{k}^\gamma -1 \leq t\leq KL +\bar{k}^\gamma$,
		\begin{align}
			\mathbb{P}\left(\varrho_{R,i}=1,\, Y_{R,i}\leq t-\epsilon^{*}_{N,KL,K^\prime L^\prime}\right)\leq \mathbb{P}\left(-\bar{k}^\gamma + \overset{\diamond}{\epsilon}\leq \tilde{M}^f_{N,i}\leq t \right)
			\leq \mathbb{P}\left(\varrho_{R,i}=1,\, Y_{R,i}\leq t+ \epsilon^{*}_{N,KL,K^\prime L^\prime}/2 \right).
		\end{align}
		Thus, by the same argument given in the proof of \cite[Theorem~5.9]{paper2}, there is a coupling between $\{\tilde{M}^f_{N,i}:\,1\leq i \leq R \}$ and $\{(\varrho_{R,i},Y_{R,i}):\, 1\leq i \leq R\}$ such that on the event $(\mathcal{E}\cup\mathcal{E}^\prime)^c$:
		\begin{align}
			&\varrho_{R,i}=1, |Y_{R,i}-\tilde{M}^f_{N,i}|\leq \epsilon^{*}_{N,KL,K^\prime L^\prime}, \quad\text{if } \tilde{M}^f_{N,i}\geq \epsilon^{*}_{N,KL,K^\prime L^\prime}\\
			&|Y_{R,i}-\tilde{M}^f_{N,i}|\leq \epsilon^{*}_{N,KL,K^\prime L^\prime}, \,\,\,\qquad\qquad\text{if }\varrho_{R,i}=1.
		\end{align}
	 As $\{Z_{R,i}\}_{1\leq i \leq R}$ and $\{S^{N,c}_v\}_{v\in V_N}$ have the same law, one can couple such that $S^{N,c}_v=Z_{R,i}$, for $v\in B_{N/KL,i}$ and $1\leq i \leq R$. Using \cite[(5.63)]{paper2}, we deduce
		\begin{align}\label{equation:4.65}
			\limsup\limits_{(N,L,K,L^\prime, K^\prime)\Rightarrow \infty} \mathbb{P}\left(\varrho_{R,\tilde{\tau}_i}=1\right)=1,
		\end{align}
		and thereby exclude that the local maximum is achieved in a box $T^{(KL)}_j$ when at the same time $\varrho_{R,j}=0$.
		Thus, there are couplings, such that outside an event of vanishing probability as $(N,L,K,L^\prime, K^\prime)\Rightarrow \infty$, we have
		\begin{align}
			\left( (\max_{v\in A^\prime_1}S^N_v -m_N)- G^{(1)}_{L,K,L^\prime,K^\prime}, \dotsc, (\max_{v\in A^\prime_p}S^N_v-m_N) -G^{(p)}_{L,K,L^\prime, K^\prime} \right)_{\infty} 
			\leq 2 \epsilon^{*}_{N,KL,K^\prime L^\prime},
		\end{align}
		which proves \autoref{theorem:coupling_proof}.
	\end{proof}
\end{thm}
Next, we prove \autoref{theorem:conv_multiple_maxima}.
\begin{proof}[Proof of \autoref{theorem:conv_multiple_maxima}:]
	By \autoref{lemma:reduction_multiple_to_3field}, \eqref{equation:4.49} and \autoref{theorem:coupling_proof}, we can reduce the proof to proving convergence of the laws $\bar{\mu}_{L,K,L^\prime,K^\prime, \underline{A}}$. Recall that we write $R=KL$. In the following, we construct random variables $\{D_{KL}(A_i): \, 1\leq i\leq p\}_{K,L\geq 0}$ that are measurable with respect to $\mathcal{F}^C\coloneqq \sigma\left(Z_{R,i}\right)_{i=1}^{R}$, so that for any $x_1,\dotsc,x_p \in \mathbb{R}$, the following limit exists
	\begin{align}\label{equation:limits_1}
		\lim\limits_{(L,K,L^\prime, K^\prime)\Rightarrow \infty}\frac{\bar{\mu}_{L,K,L^\prime,K^\prime,\underline{A}}((-\infty,x_1],\dotsc, (-\infty,x_p])}{\mathbb{E}\left[\exp(-\beta^{*}_{K^\prime,L^\prime}\sum_{i=1}^{p}D_{KL}(A_i)e^{-2 x_i})\right]
		},
	\end{align}
	and is equal to one.
	Regarding \eqref{equation:4.65}, assume $\varrho_{R,\tilde{\tau}_i}$, for $1\leq i\leq p$.
	Conditioning on $\mathcal{F}^c$, we have, for any $x_1,\dotsc, x_p\in \mathbb{R}$ ,
	\begin{align}\label{equation:limits_2}
		&\bar{\mu}_{L,K,L^\prime,K^\prime}((-\infty,x_1],\dotsc, (-\infty,x_p])= \mathbb{P}\left(G^{(i)}_{L,K,L^\prime,K^\prime }\leq x_i:\, i=1,\dotsc,p \right)\nonumber\\
		&\qquad= \mathbb{E}\left[ \prod_{i=1}^{p}\left(1- \mathbb{P}\left( \varrho_{R,j}(Y_{R,j}+2\log(KL)(\sigma_1^2-1))>x_i +2\log(KL)-Z_{R,j} | \mathcal{F}^c\right) \right)^{|T^{(KL)_i |}}  \right].
	\end{align}
	A union bound on $\mathcal{D}^c=\{\min_{1\leq i \leq R} 2\log(KL) - Z_{R,i}\geq 0\}^c$, shows that
	\begin{align}
		\limsup\limits_{(L,K)\implies \infty} \mathbb{P}(\mathcal{D})=1.
	\end{align} 
	Thus, on the event $\mathcal{D}$, and by \eqref{equation:4.50}, \eqref{equation:4.51} and \eqref{equation:D_KL}, one deduces
	\begin{align}\label{euation:prob_cond_on_coarse_field}
		\mathbb{P}\left(\varrho_{R,j} Y_{R,j}\geq 2\log(KL ) \sigma^2(0)-Z_{R,j} + x_i | \mathcal{F}^c\right) = \beta^{*}_{K^\prime, L^\prime}e^{-2 (2(1+\sigma^2(0))\log(KL) - Z_{R,j}+x_i)}
	\end{align}
	In particular, note that \eqref{euation:prob_cond_on_coarse_field} tends to zero as $KL\rightarrow \infty$. Using $e^{-\frac{x}{1-x}}\leq 1-x \leq e^{-x}$, for $x<1$, and inserting for $x$ the probability in \eqref{euation:prob_cond_on_coarse_field} with $K,L$ large, implies that there is non-negative sequence $\{\epsilon_{K,L}\}_{K,L\geq 0}$, with $\limsup\limits_{(K,L)\Rightarrow \infty} \epsilon_{K,L}=0$, such that
	\begin{align}\label{equation:4.72}
		\exp\left(-(1+\epsilon_{K,L})\beta^{*}_{K^\prime,L^\prime}e^{-2((1+\sigma^2(0))\log(KL)-Z_{R,j}+x_i)}\right)
		 &\leq \mathbb{P}\left(\varrho_{R,j} Y_{R,j}\leq 2\log(KL)\sigma^2(0)-Z_{R,j}+x_i |\mathcal{F}^c\right)\nonumber\\
		&\leq \exp\left(-(1-\epsilon_{K,L})\beta^{*}_{K^\prime,L^\prime}e^{-2((1+\sigma^2(0))\log(KL)-Z_{R,j}+x_i)}\right).
	\end{align}
	Plugging \eqref{equation:4.72} into \eqref{equation:limits_2} gives \eqref{equation:limits_1}, with
	\begin{align}\label{equation:D_KL}
			D_{K,L}(A_i)=\sum_{j\in T^{(KL)}_i}e^{-2(2(1+\sigma^2(0))\log(KL)-Z_{R,j})}.
	\end{align}
	\eqref{equation:limits_1} combined with \autoref{theorem:coupling_proof}, implies that there is a constant $\beta^{*}>0$, such that
	\begin{align}\label{equation:beta}
		\limsup\limits_{(K^\prime,L^\prime)\Rightarrow \infty}|\beta^{*}_{K^\prime,L^\prime}-\beta^{*}|=0.
	\end{align}
	Inserting \eqref{equation:beta} into \eqref{equation:limits_1}, we obtain
	\begin{align}\label{equation:4.74}
			\lim\limits_{(L,K,L^\prime, K^\prime)\Rightarrow \infty}\frac{\bar{\mu}_{L,K,L^\prime,K^\prime,\underline{A}}((-\infty,x_1],\dotsc, (-\infty,x_p])}{\mathbb{E}\left[\exp(-\beta^{*}\sum_{i=1}^{p}D_{KL}(A_i)e^{-2 x_i})\right]
		}=1.
	\end{align}
	\autoref{theorem:coupling_proof} in combination with \eqref{equation:4.74}, implies that $\{D_{KL}(A_i):\, 1\leq i \leq p \}$ converge weakly to random variables $\{D(A_i):\, 1\leq i \leq p\}$, as $K,L\rightarrow \infty$. Moreover, as the sequence of laws, $\{\bar{\mu}_{L,K,L^\prime,K^\prime,\underline{A}} \}_{L,K,L^\prime,K^\prime \geq 0},$ is tight, it follows that almost surely, $D(A_i)>0$, for $i\in \{1,\dotsc,p\}$. This concludes the proof.
\end{proof}
	\section{Appendix}
	\subsection{Gaussian comparison}
	We need a vector version of Kahane's theorem. %(for reference look J.-P. Kahane (1985). Sur le chaos multiplicatif.), i.e.
\begin{thm}\label{theorem:vector_kahane}
	Let $f\in C^2(\mathbb{R}^n;\mathbb{R}^k)$ with sub-Gaussian growth in every component of the second derivatives. Further let $\{X_i\}_{1\leq i \leq n}, \,\{Y_i\}_{1\leq i \leq n}$ be two centred Gaussian fields satisfying
	\begin{align}
		\mathbb{E}\left[Y_i Y_j \right]> \mathbb{E}\left[X_i X_j \right] \implies \frac{\partial f}{\partial x_i \partial x_j }(x) \geq 0,\quad x\in \mathbb{R},
	\end{align}
	where the inequality is to be understood component-wise.
	Then,
	\begin{align}
		\mathbb{E}\left[f(Y)\right]\leq \mathbb{E}\left[f(X)\right],
	\end{align}
	again to be understood as an inequality valid in each component.
	\begin{proof}
		The proof is an immediate adaptation of the original proof, as each component of $f$ is a function $f_i \in C^2(\mathbb{R}^n)$ with sub-Gaussian growth in its second derivatives, for which Kahane's theorem holds. In particular, each component of the map $f$ can be treated separately.
	\end{proof}
\end{thm}
This allows us to deduce a vector version of Slepian's inequality.
\begin{thm}\label{theorem:vector_Slepian}
	Let $T$ be a countable index set, $\{X_i\}_{i\in T},\, \{Y_i\}_{i\in T}$ be two centred Gaussian fields satisfying
	\begin{align}
		\mathrm{Var}\left[X_i\right]=\mathrm{Var}\left[Y_i\right]\, \forall i\in T \quad \text{and}\quad \mathbb{E}\left[X_i X_j \right]\leq \mathbb{E}\left[Y_i Y_j\right],\, \forall i,j\in T.
	\end{align}
	Then, for any disjoint collection of subsets $T_1,\dotsc, T_k \subset T$ and real numbers $x_1,\dotsc,x_k \in \mathbb{R}$,
	\begin{align}
		\mathbb{P}\left(\max_{i\in T_1}Y_i\leq x_1,\dotsc, \max_{i \in T_k} Y_i\leq x_k  \right) \leq \mathbb{P}\left(\max_{i\in T_1}X_i\leq x_1,\dotsc, \max_{i\in T_m} X_i\leq x_k \right).
	\end{align}
	\begin{proof}
		The proof is basically only a vector version of the original, which is why we just give a sketch. Assume for simplicity $|T|=n$.
		One takes a sequence of maps $f_l:\mathbb{R}^{n} \rightarrow \mathbb{R}^k$ of the form 
		\begin{align}
			f_l=\begin{pmatrix}
			\prod_{i\in  A_1} g^l_i(x_i)\\
			\prod_{i\in  A_2} g^l_i(x_i)\\
			\vdots\\
			\prod_{i\in  A_k} g^l_i(x_i)
			\end{pmatrix}
		\end{align} where $g^l_i(x_j)$  are smooth, non-increasing and converge from above  
		to $\mathbbm{1}_{(-\infty,x_j]}$. One notices that the requirements of \autoref{theorem:vector_kahane} are met, and an application of it finishes the proof.
	\end{proof}
\end{thm}

	\bibliography{literature}
	\bibliographystyle{abbrv}
	
	%\clearpage
\end{document}